\documentclass[preprint,12pt]{elsarticle}

\usepackage{amsmath}
\usepackage[utf8]{inputenc} 
\usepackage[T1]{fontenc}    
\usepackage{hyperref}       
\usepackage{url}            
\usepackage{booktabs}       
\usepackage{amsfonts}       
\usepackage{nicefrac}       
\usepackage{microtype}      
\usepackage{xcolor}         
\usepackage{multirow}
\usepackage{todonotes}
\usepackage{eurosym}

\usepackage{multicol}
\usepackage{etoolbox}
\usepackage{nomencl}
\usepackage{color, soul}

\journal{Elsevier}

\begin{document}

\begin{frontmatter}



\title{Risk-aware stochastic scheduling of multi-market energy storage systems}

\author[label1,label2]{Gabriel D. Patrón}
\author[label3]{Di Zhang}
\author[label1,label2]{Lavinia M.P. Ghilardi}
\author[label3]{Evelin Blom}
\author[label3]{Maldon Goodridge}
\author[label3]{Erik Solis}
\author[label3]{Hamidreza Jahangir}
\author[label3]{Jorge Angarita}
\author[label3]{Nandhini Ganesan}
\author[label3]{Kevin West}
\author[label2]{Nilay Shah}
\author[label1,label2]{Calvin Tsay}
 \affiliation[label1]{organization={Department of Computing, Imperial College London},
             city={London},
             postcode={SW7 2AZ},
             country={United Kingdom}}
 \affiliation[label2]{organization={Centre for Process Systems Engineering, Imperial College London},
             city={London},
             postcode={SW7 2AZ},
             country={United Kingdom}}
 \affiliation[label3]{organization={BP International Ltd},
             city={Sunbury-on-Thames},
             postcode={TW16 7BP},
             country={United Kingdom}}
             
\begin{abstract}

Energy storage promotes the integration of renewables by operating with charge and discharge policies that balance an intermittent power supply. A key challenge in this emerging sector is how to optimize the operation of storage assets given future price uncertainties and the need to recover the costs of project finance while ensuring an attractive return on equity and hedging against downside risk. This study investigates the scheduling of energy storage assets under price uncertainty, with a focus on electricity markets. A two-stage stochastic risk-constrained approach is employed, whereby electricity price trajectories or specific power markets are observed, allowing for recourse in the schedule. Conditional value-at-risk is used to quantify risk in the optimization problems; this allows for explicit specification of a probabilistic risk limit. The proposed approach is tested in an integrated  hydrogen system (IHS) and a battery energy storage system (BESS). In the joint design and operation context for the IHS, the risk constraint results in large installed unit capacities, increasing capital cost but enabling more inventory to buffer price uncertainty. In both case studies, there is an operational trade-off between risk and expected reward; this is reflected in higher expected costs (or lower expected profits) with increasing risk aversion. Despite the decrease in expected reward (up to $\approx 500 \$k$), both systems exhibit substantial benefits of increasing risk aversion (up to $\approx 1.5 \$mn$) with respect to risk-neutral settings. This work provides a general method to address uncertainties in energy storage scheduling, allowing operators to input their level of risk tolerance on asset decisions.

\end{abstract}

\begin{graphicalabstract}
\begin{figure}[t]
\centering
\includegraphics[width=1\textwidth]{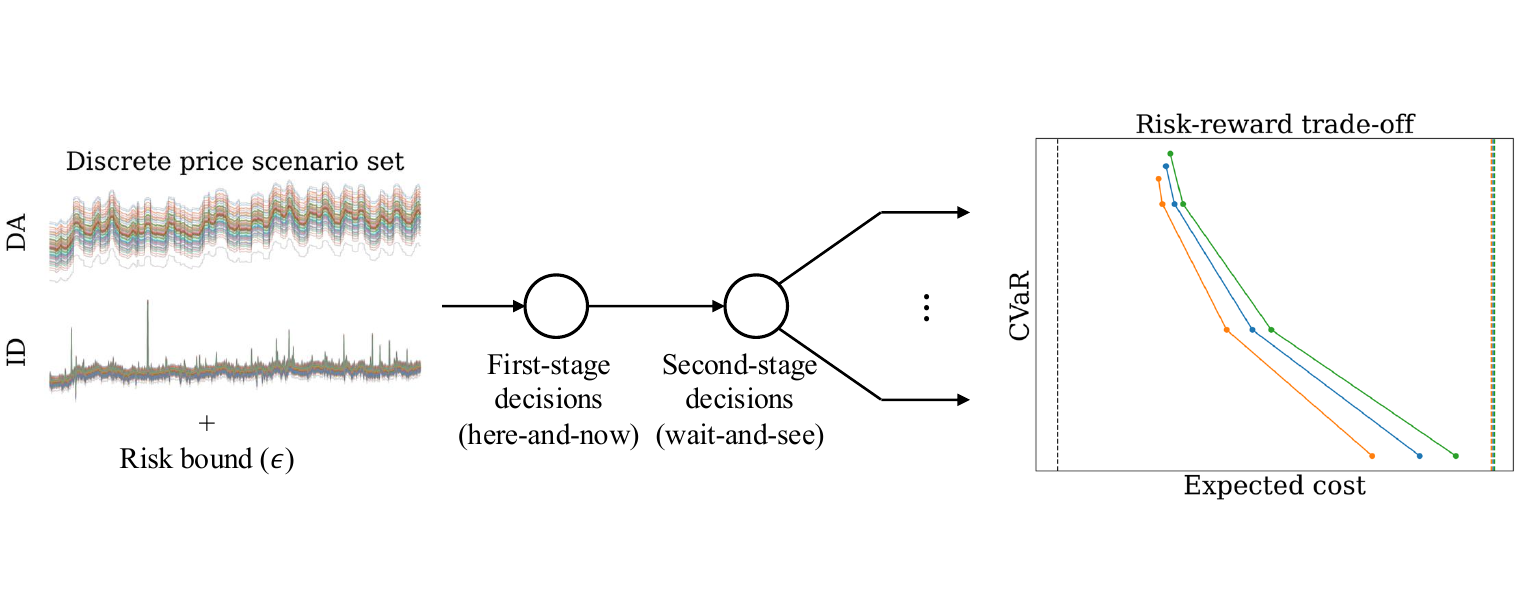}
\end{figure}
\end{graphicalabstract}

\begin{highlights}
\item Two-stage stochastic optimization is proposed for scheduling energy storage systems
\item Coordinated participation across multiple power markets is modelled
\item Unit capacities and charge/discharge schedules are jointly optimized
\item Conditional value-at-risk (CVaR) was used to limit risk owed to price uncertainty 
\item The approach was tested in integrated hydrogen and battery energy storage systems
\end{highlights}


\end{frontmatter}



\section{Introduction}
\label{sec1}

Worldwide renewable energy capacity is rapidly expanding, accounting for $38\%$ of the increase in global energy supply in 2024 \citep{iea2025}. For this transition to continue, energy storage systems will be key in abating intermittencies \citep{mallapragada2020long}, which limit the uptake of renewable generation and necessitate polluting base load energy sources such as combined-cycle power plants. Among the storage methods that have been proposed \citep{koohi2020review}, hydrogen and battery energy storage systems are prominent in various national clean energy action plans (e.g., EU~\citep{eu2023}, UK \citep{uk2025}).

Hydrogen produced using an electrolytic cell powered by renewable energy is referred to as ``green hydrogen'' owing to its low carbon intensity \citep{ajanovic2022economics}. When supplied with electricity, an electrolyzer converts water to oxygen gas in the anode and hydrogen gas in the cathode. Hydrogen is a versatile energy carrier, with a variety of uses, and a potential to further integrate variable renewable electricity into the grid \citep{mallapragada2023decarbonization}. In grid-connected settings, hydrogen can provide both short- and medium-term pressurized storage and be integrated with a fuel cell or turbine for conversion back into electricity. Excess hydrogen may also be used as a drop-in fuel~\citep{ueckerdt2021potential} for transport applications or as a feed stock to ammonia production~\citep{lee2022pathways}. Note that, while pathways toward green grid-connected hydrogen have been investigated~\citep{ricks2023minimizing}, the carbon intensity of grid-connected hydrogen production depends on the current associated generation mix. Furthermore, \mbox{\citet{babay2025forecasting}} investigate forecasting techniques for green hydrogen production under renewable variability, highlighting the need for uncertainty-aware decision support in this type of system.

In contrast to hydrogen storage, which has seen little uptake given its relatively high capital cost, battery energy storage system (BESS) processes have seen a larger installed capacity to date \citep{emrani2024comprehensive}. A BESS typically employs the electrolytic reduction of lithium ions to store electrical energy in a cell. This enables short-term storage enabled by fast charge/discharge cycles, which can satisfy electrical grid demand. Increasingly, these battery assets are being used for energy arbitrage between various electrical markets \citep{pusceddu2021synergies,krishnamurthy2017energy,nezamabadi2020arbitrage}, where price spreads are exploited to generate a profit. By responding to changing electricity prices, a BESS operator, or indeed any rapid-response storage medium operator, can capitalize on these arbitrage opportunities.

Demand response refers to the dynamic adjustment of power consumption by an electricity consumer according to the energy prices, which are time-varying (often called ``indirect'' demand response). To this end, optimal demand response scheduling has been applied in many system contexts including: residential/industrial buildings \citep{nan2018optimal,dowling2017multi}, water distribution systems \citep{oikonomou2018optimal}, and industrial air separation units \citep{tsay2019optimal}. Energy storage further allows for effective demand response by purchasing and storing energy during periods of low prices, thereby allowing for grid curtailment and consumption of stored energy during periods of high prices \citep{li2023optimal,zhang2020two,tang2019model}. Most of these approaches follow the so-called ``price-taker'' approach, relying on accurate forecasting of the markets involved as problem inputs, which are subject to uncertainties in practice. These uncertainties induce suboptimal scheduling of the demand response, which results in economic losses. The reader is referred to \citet{silva2022demand} for a review of uncertainty in demand response.

To abate the effect of these price uncertainties, the literature has turned to stochastic optimization \citep{powell2019unified,li2021review}, which allows for a distribution over electricity price estimates to be embedded into the optimal demand response scheduling problem \citep{silva2022demand}. In particular, two-stage stochastic optimization, with several readily-deployable software packages \cite{torres2022review}, provides an attractive framework for demand response decisions to be segmented into here-and-now and wait-and-see decisions. The former are implemented immediately, while the latter can be made at later time when uncertainties are realized. While a choice of single here-and-now decisions is actioned, a conditional distribution of wait-and-see decisions is obtained such that the demand response schedule is adjusted according to the realized prices. For instance, here-and-now decisions may correspond to capital installations that must be decided before electricity prices are known, while wait-and-see decisions may involve how those processes are scheduled. Despite the success of stochastic demand response, these methods optimize the wait-and-see decisions based on their expectation; this can potentially lead to poor outcomes if the price scenarios that correspond to the tails of the distribution are realized in practice. Accordingly, a risk measure such as conditional value-at-risk (CVaR) can be modelled \citep{filippi2020conditional}, optimized as an objective (penalty), or constrained.

In the context of energy storage, the joint explicit optimization of expected cost and tails risk for energy storage systems has been explored in microgrids with BESS \citep{do2021decision, herding2024risk} and natural gas storage \citep{moradi2022risk} settings for single power markets. Furthermore, this has been extended to BESS processes operating in multiple power markets \citep{nezamabadi2020arbitrage}. The joint optimization of cost expectation and risk, however, results in a multi-objective optimization problem, which is subject to an operator-defined weighting between the two objectives. The heuristic nature of the objective weights results in a lack of probabilistic guarantees on the tail risk. Notably, \citet{herding2024risk} explored the use of CVaR constraints in single-market BESS to explicitly limit the tail risk as specified in \citet{haimes1971bicriterion}; this provides a more attractive proposition to system operators and arbitrageurs that desire precise control over their exposure.

Previous works have studied risk management in energy storage for narrowly defined markets, systems, and problem conditions; to the best of the authors' knowledge, there is no work that systematically explores the effect of risk across multiple settings. Based on the literature reviewed above, there is a gap for a general-purpose risk-aware stochastic optimization scheme that can be applied across energy storage media, problem formulations, and electricity markets.  In this work, we propose a CVaR-constrained stochastic optimization framework for the scheduling of energy storage assets under electricity price uncertainty and show its applicability to multiple systems and decision structures. The proposed approach is demonstrated using two prototypical energy storage systems: integrated hydrogen system (IHS) and a BESS. Both systems are considered to participate in day-ahead (DA) and intraday (ID) power markets. While these selected storage systems are similar in their ability to charge and discharge from the grid, they provide distinct perspectives on the risk-constrained scheduling problem. The IHS system is subject to a cost minimization problem that features capital decisions (i.e., it is a joint design and operation problem). Conversely, fixed capacities are assumed for the BESS system in a profit maximization (arbitrage) context; the fixed capacities enable scheduling to be applied in both long-term scheduling and rolling horizon (i.e., feedback control) settings. For BESS, the former setting can be used determine the optimal dispatch for a flexible power purchasing agreement, while the latter setting can be used in live energy trading. The risk-reward trade-off, proportion of market participation, and choice of recourse variables are explored in both systems. The key contributions of this study are:

\begin{itemize}
    \item A CVaR-constrained two-stage stochastic optimization framework is proposed for storage assets that explicitly accounts for downside risk arising from intertemporal state coupling (e.g., in state-of-charge or hydrogen inventory).
    \item Coordinated participation of storage in coupled day-ahead and intraday markets is modelled, showing how commitments and recourse decisions interact with physical storage constraints.
    \item The framework is deployed on an integrated hydrogen system for joint design and operation under price risk as well as a battery energy storage system used for multi-market arbitrage under downside risk constraints.
    \item The approach is shown to support both long-term planning and rolling-horizon operation, including quantitative evaluation of flexible power purchasing agreement structures through delayed price observation.
\end{itemize}

The remainder of this work is organized as follows: \autoref{sec2} details stochastic optimization, CVaR modelling, and the proposed risk-constrained approach; \autoref{sec3} and \autoref{sec4} present the IHS and BESS case studies, respectively; \autoref{sec5} applies the proposed approach to the previously outlined case studies and provides analysis on its outcomes; \autoref{sec6} outlines the key takeaways and future work directions for this field.

\section{Methodology}
\label{sec2}

Deterministic opimization problems often require inputs that are not precisely known \textit{a priori} (i.e., they are uncertain). The solution to a deterministic programming problem can be suboptimal with respect to the real-life systems subject to these uncertainties, which can result in economic losses. Uncertain inputs are often forecasted at the time of optimization ($t_0$) and revealed later (e.g., once a market auction occurs at $t_\mathrm{obs}$). Stochastic programming \citep{powell2019unified} can be used to abate the effects of uncertainty and limit potential losses. In this section, we propose the use of a stochastic programming framework that can be used for the short- and long-term optimization of energy storage systems. Further, we present a formulation for constraining CVaR to probabilistically limit tail risk.

\subsection{Two-stage stochastic optimization}
\label{subsec1}

\begin{figure}[t]
\centering
\includegraphics[width=1\textwidth]{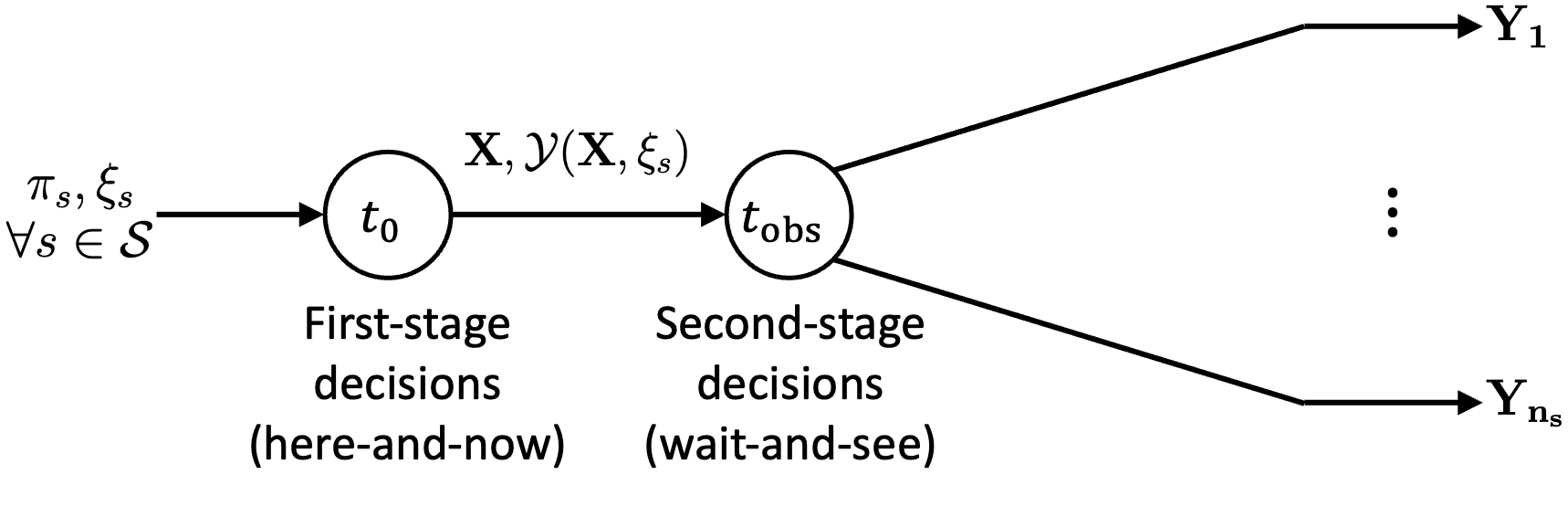}
\caption{Tree diagram for two-stage stochastic program.}\label{fig1}
\end{figure}

Stochastic programs for energy systems optimization subject to uncertain time-dependent signals (e.g., electricity market prices) are often formulated as multi-stage problems \citep{wang2022multi,barbar2022decision}, where decisions are segmented by time instances (i.e., stages) that are dependent on the solution of previous stages and values of uncertain inputs, which are dynamically revealed. These multi-stage problems require significant computational effort to solve (e.g., \citep{al2021two}), thus we deploy a two-stage approximation of the stochastic program for the optimization of energy storage systems. In the two-stage paradigm, decisions are only made at two time instances, the first-stage (`here-and-now') action and the second-stage (`wait-and-see') recourse action. While the two-stage formulation is less expressive in the temporal domain and may lead to suboptimalities, it allows for quicker solution time, ease of interpretability, and does not require successive conditional probability distributions, which are often unavailable. 

This work focuses on electricity price uncertainty in the DA and ID markets. Although the DA market clears prior to physical delivery, modelling recourse remains essential. Herein, recourse does not imply revision cleared DA commitments \textit{ex post}; rather, it represents the operator's ability to adapt subsequent operational decisions (e.g., ID participation and participation in future DA settlement periods) after additional price information and dispatch outcomes are observed over time. These dispatch decisions determine the system state entering later periods (e.g., storage inventory and available flexibility), which in turn affects the feasibility and profitability of future actions. The proposed two-stage formulation captures this intertemporal coupling, which allows risk exposure arising from commitments to be explicitly managed, thereby allowing flexibility in later market participation to hedge against adverse price realizations.

A broad formulation allows us to adapt the two-stage demand response scheduling problem to various problem settings (e.g., hydrogen or BESS) and decision structures (i.e., first- and second-stage variables). The general stochastic optimization problem is posed as follows:

\begin{equation}
\label{eq1}
\min_{\mathbf{X}\in \mathcal{X},\mathbf{Y}\in\mathcal{Y(\mathbf{X}, \mathbf{\xi})}} \mathbb{E}_{\mathbf{\xi}} \left[\mathcal{L}(\mathbf{X},\mathbf{Y},\mathbf{\xi})\right]  =\min_{\mathbf{X}\in \mathcal{X},\mathbf{Y}\in\mathcal{Y(\mathbf{X}, \mathbf{\xi})}} \mathbf{c}^{\top}\mathbf{X}+\mathbb{E}_{\mathbf{\xi}} \left[V(\mathbf{X},\mathbf{Y},\mathbf{\xi})\right],
\end{equation}
where the expected value of the objective function  $\mathcal{L}:\mathcal{X}\times\mathcal{Y(\mathbf{X}, \mathbf{\xi})}\times\Psi\rightarrow\mathbb{R}$ distribution is conditioned on the uncertainty set $\mathbf{\xi}\in\Psi$ with distribution $\mathcal{P}$ and support $\Psi$. The objective is minimized by the first-stage $\mathbf{X}\in\mathcal{X}\subset\mathbb{R}^{X}$ and second-stage $\mathbf{Y}\in\mathcal{Y(\mathbf{X},\mathbf{\xi})}\subset\mathbb{R}^{Y}$ decisions, the latter of which are dependent on the former and the uncertainty set. Assuming a linear first-stage objective (in general, we will deal with LPs herein), the LHS of \autoref{eq1} can be decomposed into first- and second-stage objectives, where $\mathbf{c}\in\mathbb{R}^{X}$ is the first-stage cost, and only the second-stage objective $V:\mathcal{X}\times\mathcal{Y(\mathbf{X}, \mathbf{\xi})}\times\Psi\rightarrow\mathbb{R}$ is conditioned on ${\mathbf{\xi}}$. This breakdown separates the `here-and-now' decisions ($\mathbf{X}$) from the `wait-and-see' decisions ($\mathbf{Y}$).

For this problem to be solved using standard optimization solvers, \autoref{eq1} must be formulated in a closed (i.e., solvable), deterministic form. Therefore, the conditional objective function term is often discretized into a finite set $\mathcal{S}=\{1,...,n_S\}$ of $n_S$ realizations of uncertainty (i.e., scenarios)\citep{kim2014guide}. This produces a large, monolithic approximation of \autoref{eq1}, referred to as the sample average approximation (SAA):

\begin{equation}
\label{eq2}
\min_{\mathbf{X}\in \mathcal{X},\mathbf{Y}\in\mathcal{Y(\mathbf{X}, \mathbf{\xi})}} \mathbf{c}^{\top}\mathbf{X}+\sum_{s\in\mathcal{S}}\pi_{s}v(\mathbf{X},\mathbf{Y}_{s},\mathbf{\xi}_s),
\end{equation}
where $\pi_s\in\mathbb{R}$ represents the probability of scenario $s\in\mathcal{S}$ corresponding to a discrete realization of the second-stage objective $v$. A schematic of the discretized two-stage decision structure is shown in \autoref{fig1} whereby the first- and second-stage decisions are made at $t_0$ and $t_\mathrm{obs}$, respectively. The observation time $t_\mathrm{obs}$ denotes the time at which the values of uncertain variables are revealed, and thus the `wait-and-see' decisions are actioned according to the scenario that is realized in practice. 
Several alternatives to the SAA formulation have been proposed, such as by using dynamic optimization~\citep{tsay2017dynamic} or surrogate models for the second-stage objective~\citep{patel2022neur2sp}. For example, our recent work~\citep{ghilardi2025integrated} explores surrogate models for integrated design and scheduling of a hydrogen process under uncertainty. 

In this formulation, we optimize on a fixed horizon with a single observation time; this abstraction is used to obtain a tractable two-stage approximation of the underlying multi-stage market process. The $t_\mathrm{obs}$ parameter approximates the point in time at which sufficient price information becomes available to materially alter operational strategy, allowing us to study how delayed information affects risk exposure and performance. We emphasize that the time partitioning does not represent a market-clearing sequence itself, but rather the timing of information availability used to parameterize recourse in a two-stage approximation; dependencies between day-ahead commitments and subsequent intraday operation are captured through intertemporal state coupling and operational constraints. The two-stage approximation simplifies temporal risk propagation relative to a full multi-stage formulation, as only a single information update is modelled. In storage systems with strong intertemporal coupling, especially those with long-horizons and inventory coupling, intermediate recourse opportunities may alter both expected performance and risk exposure. Accordingly, the present study should be interpreted as a tractable approximation to the underlying sequential decision problem. Comparison between two-stage and multi-stage formulations have been performed for unit-commitment settings in electricity markets \mbox{\citep{mahmutougullari2019value}} and the aim of the present work is to show downside risk can be effectively managed even within a computationally tractable, two-stage framework.

\subsection{Conditional value-at-risk}
\label{subsec2}

The formulation presented in \autoref{eq1} (and approximated in \autoref{eq2}) optimizes over the expected value of the conditional objective distribution; however, tail risk is another potential optimization objective to mitigate against extreme shortfall. In other words, decision makers may be concerned with the expected performance given likely circumstances, but also with worst-case performance given a more extreme scenario. 
The conditional value-at-risk (CVaR) is a common tail risk metric used in optimization owing to its convexity and coherence (i.e., monotonicity, sub-additivity, homogeneity, translational invariance)~\citep{artzner1999coherent}. For a conditional random variable such as the second-stage loss function $V$ with a cost minimization objective as in \autoref{eq1}, the right-tail (i.e., tail cost) CVaR is expressed as:

\begin{eqnarray}
\label{eq3} \mathrm{CVaR}_{\alpha}(V) = \mathbb{E}_{\xi}[V|V\geq \zeta]\\
\label{eq4} \zeta = \inf\{V \in \mathbb{R}:\mathcal{F}_V(\xi)\geq\alpha\},
\end{eqnarray}
where $\mathcal{F}_V$ is the cumulative distribution function of $V$, and $\alpha$ is a user-specified risk percentile. The variable $\zeta$ represents the value-at-risk (VaR). Alternatives to using CVaR to mitigate risk in energy systems optimization include using other risk metrics or robust optimization-based formulations. The reader is referred to \citet{rahim2022overview} for an review of robust optimization in energy grid applications. We consider CVaR in this work because it directly quantifies tail risk, which is more relevant to storage operators than a dispersion measure like variance. Unlike VaR, CVaR also accounts for losses beyond a quantile threshold and is a coherent risk measure. Finally, CVaR admits a tractable linear reformulation under the SAA in \mbox{\autoref{eq2}}, allowing it to be directly deployed in a two-stage stochastic program. Alternative paradigms such as robust optimization protect against worst-case realizations and can lead to (overly) conservative solutions.

CVaR can be reformulated into individual realizations of uncertainty to produce an explicit closed-form optimization formulation, as shown by \citet{rockafellar2000optimization}. This reformulation yields the following:

\begin{equation}
\label{eq5}
\mathrm{CVaR}_{\alpha}(v) = \zeta + \frac{1}{1-\alpha}\sum_{s\in \mathcal{S}}\pi_{s}[v(\mathbf{X},\mathbf{Y}_{s},\mathbf{\xi}_s) -\zeta]^{+},
\end{equation}
where \autoref{eq5} sums the difference between the VaR and the objectives found in the SAA scenarios that exceed the VaR itself (i.e., using the positive part operator $[\cdot]^{+}=\mathrm{max}\{\cdot,0\}$). This expression can be further reformulated to avoid a bi-level optimization problem (i.e., including the max operator) as follows:

\begin{eqnarray}
\label{eq6} \mathrm{CVaR}_{\alpha}(v) \equiv\zeta + \frac{1}{1-\alpha}\sum_{s\in \mathcal{S}}\pi_{s}\eta_{s}\\
\label{eq7} v(\mathbf{X},\mathbf{Y}_{s},\mathbf{\xi}_s) -\zeta \leq \eta_{s} ; \quad \forall s\in\mathcal{S} \\
\label{eq8} \eta_{s} \geq 0 ;\quad \forall s\in\mathcal{S},
\end{eqnarray}
where $\eta_s$ is a non-negative auxiliary variable introduced in the reformulation.

\subsection{Risk-constrained two-stage stochastic optimization}
\label{subsec3}

Combining the SAA objective in \autoref{eq2} with the CVaR formulation in \autoref{eq6}--\autoref{eq8} yields the resulting risk-constrained two-stage problem:

\begin{eqnarray}
\label{eq9}
\min_{\mathbf{X}\in \mathcal{X},\mathbf{Y}_s\in\mathcal{Y}} \mathbf{c}^{\top}\mathbf{X}+\sum_{s\in\mathcal{S}}\pi_{s}v(\mathbf{X},\mathbf{Y}_{s},\mathbf{\xi}_s) \nonumber \\
s.t. \nonumber\\
\zeta + \frac{1}{1-\alpha}\sum_{s\in \mathcal{S}}\pi_{s}\eta_{s} \leq\epsilon \\
\mathbf{f}(\mathbf{X},\mathbf{Y}_{s},\mathbf{\xi}_s) \leq 0 &; \quad \forall s\in\mathcal{S} \nonumber\\
v(\mathbf{X},\mathbf{Y}_{s},\mathbf{\xi}_s) -\zeta \leq \eta_{s}&; \quad \forall s\in\mathcal{S} \nonumber\\
\eta_{s} \geq 0&; \quad \forall s\in\mathcal{S}, \nonumber
\end{eqnarray}
where the user-specified parameter $\epsilon$ reflects risk aversion by providing an upper bound to CVaR. In practice, this upper bound represents the user-defined maximum allowable $\alpha^{th}$ percentile tail loss, thus limiting the potential for extreme shortfall. To choose a value for $\epsilon$, one may first solve the risk-neutral stochastic problem to determine the nominal CVaR; the operator can then choose to impose \mbox{$\epsilon$} value lower than the nominal CVaR to reflect the aggressiveness of their energy trading strategy (e.g., through successive tightening of the constraint). A trade-off exists between expected cost and CVaR, where a larger bound results in lower expected costs (i.e., a less conservative formulation introduces more risk, but a potentially higher  expected reward). A risk-aware formulation where the CVaR expression is embedded into the objective function with a weighting factor is also used in the literature \citep{do2021decision, moradi2022risk, herding2024risk,nezamabadi2020arbitrage}; however, this may be less useful from an operator's perspective, as it does not directly enable specifying an upper risk limit. The risk-constrained approach presented in \autoref{eq9} provides an explicit risk bound, which provides a probabilistic guarantee of limited shortfall. The general constraints $\mathbf{f}:\mathcal{X}\times\mathcal{Y(\mathbf{X},\mathbf{\xi})}\times\Psi\rightarrow\mathbb{R}^{f}$ correspond to the system model, which further imposes operational constraints on the optimization problem and is instantiated in the scenario set $\mathcal{S}$. 

The formulation in \autoref{eq9} is general and can be deployed in a variety of energy storage systems, operational settings, and  uncertain inputs as will be shown in the forthcoming case studies. In this work, we  take $\mathbf{f}(\cdot)$ to be a linear (or linearized) process model, thus all optimization problems are linear programs (LPs).

\section{Integrated hydrogen system}
\label{sec3}

\begin{figure}[t]
\centering
\includegraphics[width=1\textwidth]{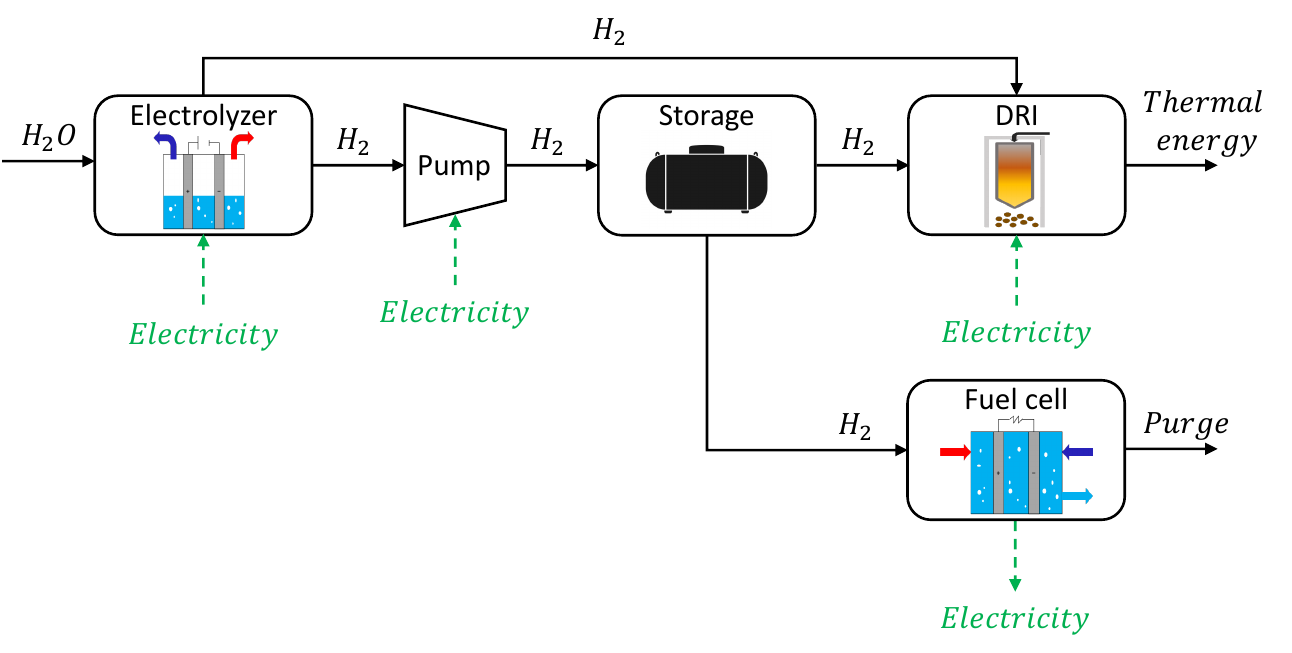}
\caption{Schematic of the integrated hydrogen system.}\label{fig2}
\end{figure}

In this section, we adapt an integrated hydrogen system (IHS) case study based on \citet{tsay2023integrating}. This system, shown in \autoref{fig2}, comprises an electrolyzer that generates hydrogen to satisfy the fixed demand of a direct reduced iron (DRI) furnace. Alternatively, the hydrogen can be compressed for storage, where it can later deployed in the DRI furnace or converted back into electricity via a fuel cell. There is no fixed electrical load for this process; the electrolyzer, pump, and DRI can purchase variable amounts of electricity provided that there is enough new and stored hydrogen to satisfy the DRI demand. The electricity prices are unknown at the time the plant is built and are revealed when the IHS begins operation or some time thereafter. Accordingly, recourse allows for capital decisions to be robust to future unknown electricity prices at construction time. The IHS model from \mbox{\citet{tsay2023integrating}} admits a linear program for the stochastic optimization framework; accordingly, several thermodynamic effects (e.g., efficiencies, heat integration, pressure-temperature coupling) are omitted or approximated. These assumptions may affect the quantitative valuation of hydrogen storage flexibility, and capacity decisions should be treated as directional insights within a tractable risk-aware framework rather than exact plant-specific prescriptions. While we take the nominal IHS model for the present study; we refer the reader to \mbox{\citet{tsay2023integrating}} for a sensitivity analysis on both key physical and capital cost parameters.

The set of units in this system is denoted $\mathcal{U}=\{\mathrm{elec},\mathrm{stor},\mathrm{heat},\mathrm{comp},fc\}$. A fixed one-year operating horizon is assumed, which is discretized into hour-long intervals $\mathcal{T}=\{t_0,...,t_\mathrm{obs}, ..., t_f\}$ where $t_f=8760$  hours. Further, we assume participation in the day-ahead and intraday power markets, denoted $\mathcal{M}=\{DA,ID\}$. We assume that the true price signal is unknown at optimization time $t_0$, but  estimates/predicted scenarios are available, motivating a stochastic optimization approach. Although the multi-stage setting is the  most accurate representation of the energy market (i.e., prices are revealed continually through time), we approximate this problem using two stages, where the true price trajectories are revealed at a single point in time. Accordingly, once sufficient time $t_\mathrm{obs}$ has elapsed, the true price trajectory is observed, allowing for recourse actions to be taken. We define the partitioned time periods as $\mathcal{T}_0 = \{t_0,...,t_\mathrm{obs}-1\}$ and $\mathcal{T}_\mathrm{obs} = \mathcal{T} \setminus\mathcal{T}_0$ (i.e., $\mathcal{T}=\mathcal{T}_0 \cup \mathcal{T}_\mathrm{obs}$).

The risk-constrained two-stage problem for this system has the first-stage decisions $\mathbf{X} = \begin{bmatrix} \mathbf{C} & \mathbf{d_{t_{0}}} \end{bmatrix}^{\top}$, where $\mathbf{C} \in\mathcal{C}\subset\mathbb{R}^{|\mathcal{U}|}$ are the unit capacities and $\mathbf{d_{t_{0}}}\in\mathcal{D}_{t_{0}}\subset\mathbb{R}^{|\mathcal{U}\times \mathcal{T}_0 \times \mathcal{M}| }$ are the charge-discharge decisions from/to the grid before the price trajectory can be observed; positive and negative domains represent charging and discharging, respectively, as the default optimization convention in \autoref{eq9} is cost minimization. The second stage decisions are the charge/discharge dispatch decision after price trajectory observation $\mathbf{Y}=\mathbf{d_{t_\mathrm{obs}}}\in\mathcal{D}_{t_\mathrm{obs}}(\mathbf{C}, \mathbf{d_{t_{0}}})\subset\mathbb{R}^{|\mathcal{U}\times \mathcal{T}_\mathrm{obs} \times \mathcal{M}| }$. The first-stage price vector is concatenated as $\mathbf{c}\in\mathbb{R}^{ |\mathcal{U}|+|\mathcal{U} \times \mathcal{T}_0 \times \mathcal{M}|}$ where they are partitioned into capital and operating prices $\mathbf{c} = \begin{bmatrix} \mathbf{P_{cap}} & \mathbf{P_{t_{0}}} \end{bmatrix}^{\top}$. The second-stage uncertain price vector is the multivariate distribution $\mathbf{\xi}=\mathbf{P_\mathrm{obs}}\in\Psi\subset\mathbb{R}^{|\mathcal{U} \times \mathcal{T}_\mathrm{obs} \times \mathcal{M}|}$. 

The IHS is modelled as a price-taking participant in both the DA and ID electricity markets. At any optimization instance, the operator schedules electricity purchases and sales for the electrolyzer, heater, compressor, and fuel cell operation in both markets. In the two-stage formulation, all DA and ID dispatch decision over the initial horizon $\mathcal{T}_{0}$ (i.e., prior to observation time) are treated as here-and-now, whereas the dispatch decisions over the remaining horizon $\mathcal{T}_\mathrm{obs}$ are modelled as wait-and-see (recourse) decisions conditioned on the realized price trajectories and system states (i.e., hydrogen inventory and available flexibility).

The model $\mathbf{f}$ describing the system is outlined next; a detailed description of the model development can be found in \citet{tsay2023integrating} along with a sensitivity analysis on key parameters. Model parameters are given in \autoref{tab1}.

\begin{table}[t]
\centering
\begin{tabular}{lllll}
\toprule
& Parameter & Symbol & Value & Units\\ 
\midrule
Electrolyzer & Rectifier efficiency & $L_{AC/DC} $ & 1.05 & - \\ 
& Auxiliary power consumption & $L_\mathrm{aux} $ & 0.05 & - \\
& Degradation factor & $L^\mathrm{deg}_\mathrm{elec}$ & 0.9142 & - \\
& Efficiency & $L_\mathrm{elec}$ & 0.05 & MW/kg \\
& Pressure & $p_\mathrm{elec}$ & 1 & MPa \\
Storage & temperature & $T_\mathrm{stor}$ & 298 & K \\
& Lower pressure bound & $p^{lB}_\mathrm{stor}$ & 2 & MPa \\
& Upper pressure bound & $p^{UB}_\mathrm{stor}$ & 20 & MPa \\
& Compressibility factor & $Z$ & 1.07 & - \\
Heater & Specific energy & $e_p$ & 11.82 & MJ/kg \\
& Heater efficiency & $\eta_\mathrm{heat}$ & 0.75 & - \\
Compressor & Efficiency & $\eta_\mathrm{comp}$ & 0.7 & - \\
Fuel cell & Inverter efficiency & $L_{DC/AC} $ & 0.95 & - \\ 
& Auxiliary power consumption & $L_\mathrm{aux} $ & 0.05 & - \\
& Degradation factor & $L^\mathrm{deg}_\mathrm{elec}$ & 0.9142 & - \\
& Voltage & $V$ & 0.7 & V \\

\bottomrule
\end{tabular}
\caption{IHS model parameters.}\label{tab1}
\end{table}

\subsection{Electrolyzer}
\label{subsec4}

The electrolyzer is parametrized using a linear scaling law to model the rectification between grid AC and plant DC current. Fixed efficiency factors are used for the DC/AC conversion and auxiliary equipment DC power consumption:

\begin{equation}
\label{eq10}
\sum_{\forall m \in M} d_\mathrm{elec}^{AC}(t,m)= (L_{AC/DC}+L_\mathrm{aux})d_\mathrm{elec}^{DC}(t);\quad\forall t \in \mathcal{T},
\end{equation}
where $d_\mathrm{elec}^{AC}$ and $d_\mathrm{elec}^{DC}$ (MW) are, respectively, the AC and DC electrolyzer power consumptions at hour $t$ and from market $m$. Electrolyzer efficiency and degradation factors are used to model the conversion of DC electricity to hydrogen gas:

\begin{equation}
\label{eq11}
F_\mathrm{elec}^{H_{2}}(t) = \frac{L_\mathrm{elec}^\mathrm{deg}d_\mathrm{elec}^{DC}(t)}{L_\mathrm{elec}};\quad\forall t \in \mathcal{T},
\end{equation}
where $F_\mathrm{elec}^{H_{2}}$ (kg/h) denotes the hydrogen production flowrate, which is split between the storage and heating units, i.e.:

\begin{equation}
\label{eq12}
F_\mathrm{elec}^{H_{2}}(t) = F_\mathrm{stor}^\mathrm{in}(t)+F_\mathrm{heat}^\mathrm{in,elec}(t);\quad\forall t \in \mathcal{T}.
\end{equation}

A constraint is imposed such that the DC power consumption of the electrolyzer does not exceed the unit capacity $C_\mathrm{elec}$ (MW):

\begin{equation}
\label{eq13}
d_\mathrm{elec}^{DC}(t) \leq C_\mathrm{elec};\quad\forall t \in \mathcal{T}.
\end{equation}

Finally, a power rate-of-change (i.e., ramping) constraint is imposed on the DC power hourly consumption:

\begin{equation}
\label{eq14}
|d_\mathrm{elec}^{DC}(t) -d_\mathrm{elec}^{DC}(t-1)| \leq 0.2 C_\mathrm{elec};\quad\forall t \in \mathcal{T}\setminus\{t_{0}\}.
\end{equation}

\subsection{Storage}
\label{subsec5}

The storage inventory is modelled by mass balances from one time period to the next: 

\begin{equation}
\label{eq15}
I_\mathrm{stor}(t) = I_\mathrm{stor}(t-1)+(F_\mathrm{stor}^\mathrm{in}(t)-F_\mathrm{stor}^\mathrm{out}(t))\Delta t;\quad \forall t \in \mathcal{T}\setminus\{t_{0}\},
\end{equation}
where $I_\mathrm{stor}$ (kg), $F^\mathrm{in}_\mathrm{stor}$ (kg/h), and $F^\mathrm{out}_\mathrm{stor}$ (kg/h) are the time-dependent hydrogen inventory, inlet flowrate, and outlet flowrate. The outlet storage flowrate is split between the fuel cell and heating units, giving the mass balance:

\begin{equation}
\label{eq16}
F_\mathrm{stor}^\mathrm{out}(t) = F_{fc}^{H_2}(t)+F_\mathrm{heat}^\mathrm{in,stor}(t);\quad\forall t \in \mathcal{T}.
\end{equation}

The ideal gas law is used to model storage pressure by incorporating a compressibility factor calculated at the centroid of the pressure bounds and at the isothermal storage temperature in \autoref{tab1}. We further assume that sufficient capacity is required to accommodate for the upper pressure bound at the storage temperature; this is reflected in the hydrogen density $\rho$ (kg/$\mathrm{m^3}$) term:

\begin{equation}
\label{eq17}
p_\mathrm{stor}(t) = p_\mathrm{stor}^{LB} + Z\frac{I_\mathrm{stor}(t)\rho(T_\mathrm{stor},p_\mathrm{stor}^{UB})RT_\mathrm{stor}}{M_{H_2}C_\mathrm{stor}};\quad \forall t \in \mathcal{T},
\end{equation}
where $R$ (J/mol/K) is the ideal gas constant and $M_{H_2}$ (kg/mol) is the molar mass of hydrogen. The storage inventory at any given time is constrained using the following bounds:

\begin{equation}
\label{eq18}
p_\mathrm{stor}^{LB} \leq p_\mathrm{stor}(t) \leq p_\mathrm{stor}^{UB};\quad\forall t \in \mathcal{T}.
\end{equation}

\subsection{Heater}
\label{subsec6}

The heater duty $d_\mathrm{heat}$ (MW) is dependent on hydrogen throughput from both electrolyzer $F^\mathrm{in, elec}_\mathrm{heat}$ (kg/s) and storage $F^\mathrm{in, stor}_\mathrm{heat}$ (kg/s) units whereby a unit-averaged specific energy $e_p$ (MJ/kg) is used:

\begin{equation}
\label{eq19}
\sum_{\forall m \in M} d_\mathrm{heat}(t,m) = \frac{1}{\eta_\mathrm{heat}}\sum_{i \in\{\mathrm{elec,stor}\}}F_\mathrm{heat}^{\mathrm{in},i}(t)e_p ;\quad\forall t \in \mathcal{T}.
\end{equation}

A constraint is imposed on the inlet stream to the heater to ensure the constant DRI furnace hydrogen demand is met:

\begin{equation}
\label{eq20}
 F_\mathrm{heat}^\mathrm{in,elec}(t)+F_\mathrm{heat}^\mathrm{in,stor}(t)\geq 150000;\quad\forall t \in \mathcal{T}.
\end{equation}

A constraint is imposed such that the heater energy consumption does not exceed its capacity $C_\mathrm{heat}$:

\begin{equation}
\label{eq21}
\sum_{\forall m \in M} d_\mathrm{heat}(t,m) \leq C_\mathrm{heat};\quad\forall t \in \mathcal{T}.
\end{equation}

\subsection{Compressor}
\label{subsec7}

The compressor duty $d_\mathrm{comp}$ is dependent on the hydrogen flowrate being processed assuming single-stage isothermal operation; this model is linearized using a first-order Taylor series expansion around $F\mathrm{_{stor}^{in}} = 0$ at the centroid of the storage pressure bounds \citep{tsay2023integrating}, which results in the following expression: 

\begin{equation}
\label{eq22}
\sum_{\forall m \in M} d_\mathrm{comp}(t,m) = \frac{RT_\mathrm{elec}}{\eta_\mathrm{comp}}F_\mathrm{stor}^\mathrm{in}(t)ln\left( \frac{p_\mathrm{stor}^{LB} +p_\mathrm{stor}^{UB}}{2p_\mathrm{elec}} \right);\quad\forall t \in \mathcal{T}.
\end{equation}

The errors introduced by this approximation, and their effects on the optimization results, were found to be relatively small by \citep{tsay2023integrating}. 
A constraint is imposed such that the compressor duty does not exceed its capacity $C_\mathrm{comp}$:

\begin{equation}
\label{eq23}
\sum_{\forall m \in M} d_\mathrm{comp}(t,m) \leq C_\mathrm{comp};\quad\forall t \in \mathcal{T}.
\end{equation}

\subsection{Fuel cell}
\label{subsec8}

The fuel cell energy output $d^{DC}_{fc}$ (MW) is denoted as a negative value since the model convention is positive sign for charging and conversely negative for discharging. This output is approximated as a function the hydrogen feed $F_{fc}^{H_2}$ (kg/h):

\begin{equation}
\label{eq24}
-d_{fc}^{DC}(t) = 2F \frac{F_{fc}^{H_2}(t)}{M^{H_2}}VL_{fc}^\mathrm{deg};\quad\forall t \in \mathcal{T},
\end{equation}
where $F$ (C/mol) is Faraday's constant and $M^{H_2}$ (kg/mol) is the molar mass of hydrogen gas. 

Conversely to the electrolyzer, an inverter is needed to convert DC $d_{fc}^{DC}$ (MW) to AC $d_{fc}^{AC}$ (MW) power in the fuel cell as:

\begin{equation}
\label{eq25}
\sum_{\forall m \in M} d_{fc}^{AC}(t,m) = (L_{DC/AC}-L_\mathrm{aux})d_{fc}^{DC}(t);\quad\forall t \in \mathcal{T}.
\end{equation}

A constraint is imposed such that the fuel cell energy production does not exceed its capacity $C_{fc}$ (MW):

\begin{equation}
\label{eq26}
-d_{fc}^{DC}(t) \leq C_{fc};\quad\forall t \in \mathcal{T}.
\end{equation}

\subsection{Cost function}
\label{subsec9}

\begin{table}[t]
\centering
\begin{tabular}{llll}
\toprule
Unit ($i\in U$) & $P_{\mathrm{cap},i}$ & $W_i$ & $O_i$ \\ 
\midrule
Electrolyzer stack & 150 \$/kW & 0.13 & 2 \$/kW/y  \\ 
electrolyzer auxiliary & 250 \$/kW & 0.08 & - \\
Storage & 1000 \$/kg & 0.1 & 10 \$/kg/y \\
Heater & 50 \$/kW & 0.13 & 1 \$/kW/y  \\
Compressor & 50 \$/kW & 0.13 & 1 \$/kW/y  \\
Fuel cell stack & 150 \$/kW & 0.13 & 2 \$/kW/y  \\
Fuel cell auxiliary & 250 \$/kW & 0.08 & -  \\

\bottomrule
\end{tabular}
\caption{IHS capital cost parameters.}\label{tab2}
\end{table}

The cost function for the hydrogen system defines the objective function and comprises capital and operating costs:

\begin{equation}
\label{eq27}
\mathcal{L}_{IHS}(\mathbf{X}, \mathbf{Y},  \mathbf{c},\xi)= J^\mathrm{cap}_{IHS}(\mathbf{C},\mathbf{P_{cap}})+J^\mathrm{op}_{IHS}(\mathbf{d_{t_{0}}},\mathbf{d_{t_{obs}}},\mathbf{P_{t_{0}}},\mathbf{P_{t_{obs}}}).
\end{equation}

The capital cost is comprised of the unit capacity costs with a lifetime annualizing factor $\mathbf{W}\in\mathbb{R}^{|U|}$ and equipment maintenance cost $\mathbf{O}\in\mathbb{R}^{|U|}$:

\begin{equation}
\label{eq28}
J^\mathrm{cap}_{IHS}(\mathbf{C},\mathbf{P_{cap}}) = (\mathbf{P_{cap}}\odot\mathbf{W}+\mathbf{O})^{\top}\mathbf{C},
\end{equation}
where capital cost factors are given in \autoref{tab2} and $\odot$ denotes the Hadamard product. Note that \citet{ghilardi2025integrated} provides a sensitivity analysis to these cost factors as prices decrease with more mature electrolyzer and fuel cell technologies.

The operating cost is the time-accrued product of the net charge/discharge of all units and the electricity prices in day-ahead and intraday markets:

\begin{equation}
\label{eq29}
J^{op}_{IHS}(\mathbf{d_{t_{0}}},\mathbf{d_{t_{obs}}},\mathbf{P_{t_{0}}},\mathbf{P_{t_{obs}}})=\mathbf{P_{t_{0}}}^\top \mathbf{d_{t_{0}}} + \mathbf{P_{t_{obs}}}^\top \mathbf{d_{t_{obs}}}. \end{equation}

When discretized in time, the overall sample-average approximation in \autoref{eq9} for the IHS problem reduces to:

\begin{equation}
\label{eq30}
(\mathbf{P_{cap}}\odot\mathbf{W}+\mathbf{O})^{\top}\mathbf{C}+\sum_{m\in\mathcal{M}}\sum_{t\in\mathcal{T}_0}P_{m,t}d_{m,t}+\sum_{m\in\mathcal{M}}\sum_{t\in\mathcal{T}_\mathrm{obs}}\sum_{s\in\mathcal{S}}\pi_{s}P_{m,t,s}d_{m,t,s}. 
\end{equation}

\section{Battery Energy Storage System}
\label{sec4}

\begin{figure}[t]
\centering
\includegraphics[width=1\textwidth]{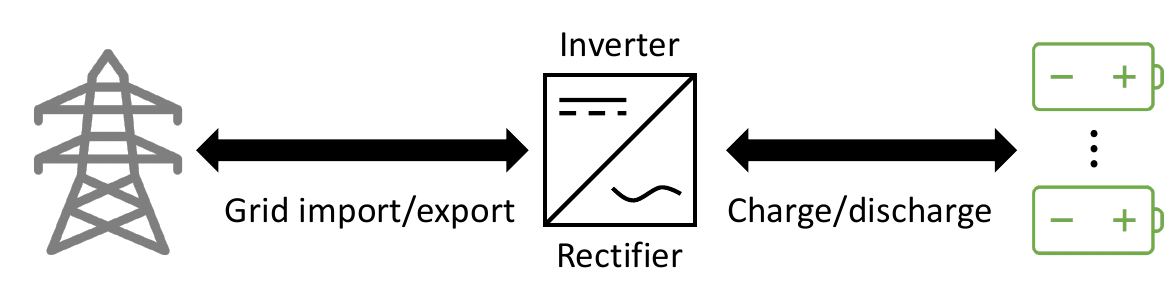}
\caption{Schematic of the battery energy storage system.}\label{fig3}
\end{figure}

This section introduces a battery energy storage system (BESS) model as a second case study to further explore the advantages of the CVaR-constrained formulation in storage settings. The system, shown in \autoref{fig3}, comprises a 50 MWh battery with fixed capacity that buys and sells electricity from the grid through an inverter/rectifier to convert current types. This battery system can be used in a variety of storage settings; however, we focus on the flexible power purchasing agreement (PPA) setting and live energy trading settings. 

As with the IHS, a fixed one-year operating horizon is discretized into hour-long intervals $\mathcal{T}=\{t_0,...,t_\mathrm{obs}, ..., t_f\}$ where $t_f=8760$  hours with participation in the day-ahead and intraday markets $\mathcal{M}=\{DA,ID\}$. The true price signal is unknown at optimization time $t_0$ and once sufficient time $t_\mathrm{obs}$ has elapsed, the true price trajectory is observed, allowing for recourse to occur. For consistency, we again define $\mathcal{T}_0 = \{t_0,...,t_\mathrm{obs}-1\}$ and $\mathcal{T}_\mathrm{obs} = \mathcal{T} \setminus\mathcal{T}_0$ (i.e., $\mathcal{T}=\mathcal{T}_0 \cup \mathcal{T}_\mathrm{obs}$).

The risk-constrained two-stage problem for this system has the first-stage decisions $\mathbf{X} = \begin{bmatrix} \mathbf{c_{t_{0}}} & \mathbf{d_{t_{0}}} \end{bmatrix}^{\top}$ where $\mathbf{c_{t_0}} \in\mathcal{C}_{t_0}\subset\mathbb{R}^{|\mathcal{T}_{0} \times \mathcal{M}|}$ are the charge and $\mathbf{d_{t_{0}}}\in\mathcal{D}_{t_{0}}\subset\mathbb{R}^{| \mathcal{T}_0 \times \mathcal{M}| }$ are the discharge decisions from the grid before the price trajectory can be observed.  The second stage decisions are the charge/discharge dispatch decision after price trajectory observation $\mathbf{Y}=\begin{bmatrix} \mathbf{c_{t_{obs}}} & \mathbf{d_{t_{obs}}} \end{bmatrix}^{\top}$ where $\mathbf{c_{t_{obs}}}\in\mathcal{C}_{t_\mathrm{obs}}(\mathbf{c_{t_{0}}}, \mathbf{d_{t_{0}}})\subset\mathbb{R}^{| \mathcal{T}_\mathrm{obs} \times \mathcal{M}| }$ and $\mathbf{d_{t_\mathrm{obs}}}\in\mathcal{D}_{t_\mathrm{obs}}(\mathbf{c_{t_{0}}}, \mathbf{d_{t_{0}}})\subset\mathbb{R}^{| \mathcal{T}_\mathrm{obs} \times \mathcal{M}| }$. For convenience and, in contrast to the IHS, both charge and discharge variables have a positive convention. The first-stage price vector is $\mathbf{c}=\mathbf{P_{t_{0}}}\in\mathbb{R}^{| \mathcal{T}_0 \times \mathcal{M}|}$ and second-stage uncertain price vector is the multivariate distribution $\mathbf{\xi}=\mathbf{P_\mathrm{obs}}\in\Psi\subset\mathbb{R}^{|\mathcal{T}_\mathrm{obs} \times \mathcal{M}|}$.

The BESS is modelled as a price-taking participant in both DA and ID electricity markets. DA commitments represent forward energy postions taken prior to price realization, while ID trading enables rebalancing and exploitation of price spreads as additional market information becomes available. In the two-stage formulation, DA and ID early dispatch decisions over the initial horizon $\mathcal{T}_{0}$ are treated as here-and-now decisions, whereas subsequent dispatch actions constitute recourse conditioned on the realized prices and system state-of-charge. This formulation can be deployed both in an open-loop planning setting and in a rolling-horizon (online) implementation. In the planning setting, the two-stage structure captures delayed price observations associated with contractual agreements such as flexible PPAs. In the rolling horizon setting, the same formulation is repeatedly re-solved as new price forecasts become available, enabling risk-constrained real-time trading while preserving intertemporal coupling through the battery state. In both cases, the CVaR constraint limits downside risk arising from aggressive DA-ID arbitrage strategies under price uncertainty. These are discussed in \autoref{subsec15} and \autoref{subsec16}, respectively, and have different mechanisms for revealing price uncertainties.

The model $\mathbf{f}$ describing the BESS is outlined next. For this case, the value of proprietary model parameters listed in \autoref{tab3} are not shared. However, similar linear models based on efficiency factors and Coulomb counting can be readily found in the literature (e.g., \citep{alavijeh2024optimal,nair2021analysis,ng2009enhanced}).

\begin{table}[t]
\centering
\begin{tabular}{lll}
\toprule
Parameter & Symbol & Units\\ 
\midrule
Battery hourly self-discharge & $\sigma_\mathrm{batt} $ & - \\ 
Round-trip efficiency & $\eta_\mathrm{batt} $ & - \\
Inverter size & $C_\mathrm{inv}$ & kW \\
Degradation constant & $\eta_\mathrm{deg}$ & - \\
Maximum daily cycles & $C_\mathrm{max}$ & - \\
\bottomrule
\end{tabular}
\caption{BESS model parameters.}\label{tab3}
\end{table}

\subsection{State of charge}
\label{subsec10}

An energy balance is used to model the state of charge $SOC$ (kWh) of the battery system:

\begin{multline}
\label{eq31}
SOC(t) = (1-\sigma_\mathrm{batt})SOC(t-1) + \eta_\mathrm{batt}   \sum_{m \in M}c(t,m)  - \\ \frac{1}{\eta_\mathrm{batt}}  \sum_{m \in M}d(t,m) ;\quad\forall t \in \mathcal{T} \setminus \{t_0\},
\end{multline}
where the charge and discharge decisions are denoted as $c$ and $d$ ($kW$), respectively. Battery self-discharge as well as round-trip efficiency are modelled using linear factors. Furthermore, the battery degrades over time, and the state of charge cannot exceed the battery capacity, as defined by the state of health $SOH$ (kWh) at a given point in time:

\begin{equation}
\label{eq32}
0 \leq SOC(t) \leq SOH(t);\quad\forall t \in \mathcal{T}. 
\end{equation}
    
\subsection{Capacity}
\label{subsec11}

Bounds are imposed on the charge and discharge decisions at every time and in every market such that they are within the limits imposed by the inverter size:

\begin{eqnarray}
\label{eq33} 0 \leq c(t,m) \leq I_{c}(t,m)C_\mathrm{inv} ;\quad\forall (t,m) \in \mathcal{T}\times \mathcal{M}\\
\label{eq34} 0 \leq d(t,m) \leq I_{d}(t,m)C_\mathrm{inv} ;\quad\forall (t,m) \in \mathcal{T}\times \mathcal{M},
\end{eqnarray}
where $I$ are binary variables to enforce no simultaneous charging and discharging in a given market and time. As separate charging and discharging variables are introduced in the BESS model, we must enforce this condition as follows:

\begin{equation}
\label{eq35}
I_{c}(t,m) + I_{d}(t,m)\leq 1 ;\quad\forall (t,m) \in \mathcal{T}\times \mathcal{M}.
\end{equation}

To ensure the total charge/discharge across all markets does not exceed the inverter/rectifier capacities, we impose the following constraints, which account for net charging or discharging:

\begin{eqnarray}
\label{eq36} \sum_{\forall m \in M}c(t,m)-d(t,m) \leq C_\mathrm{inv} ;\quad\forall t \in \mathcal{T}\\
\label{eq37} \sum_{\forall m \in M}c(t,m)-d(t,m) \leq -C_\mathrm{inv} ;\quad\forall t \in \mathcal{T}.
\end{eqnarray}

\subsection{State of health}
\label{subsec12}

With every cycle, the battery loses capacity, such that it is not able to return to its original maximum capacity; this degradation accrues over time as follows:

\begin{equation}
\label{eq38}
SOH(t) = SOH(t-1) - \eta_\mathrm{deg}\sum_{m \in M}d(t,m)   ;\quad\forall t \in \mathcal{T} \setminus \{t_0\}.
\end{equation}

To limit degradation, a cycling constraint is imposed on the battery such that a pre-specified the number of cycles performed per year is bounded. This constrains the total discharge by imposing a cumulative hourly limit on the cycles as dictated by the $SOH$:

\begin{equation}
\label{eq39}
\frac{1}{\eta_\mathrm{batt}}  \sum_{t \in \mathcal{T}} \sum_{m \in \mathcal{M}} d(t,m) \leq \frac{1}{24} C_\mathrm{max} \sum_{t \in \mathcal{T}} SOH(t).
\end{equation}

\subsection{Cost function}
\label{subsec13}

For the BESS system, we aim to maximize system profits or, to follow convention from \autoref{eq9}, minimize negative losses, i.e.:

\begin{equation}
\label{eq40}
\mathcal{L}_{BESS}(\mathbf{X,Y,c},\xi)=-(\mathbf{P_{t_{0}}}^\top (\mathbf{c_{t_{0}}}-\mathbf{d_{t_{0}}}) + \mathbf{P_{t_{obs}}}^\top (\mathbf{c_{t_{obs}}}-\mathbf{d_{t_{obs}}})).
\end{equation}

This system does not contain a capital cost term, as we found capital cost problems to be unbounded toward maximizing battery capacity for this linear formulation in the absence of market feedback. Therefore, the capacity is fixed according to available resources, and we only consider the optimal stochastic scheduling problem.

When discretized in time, the overall sample-average approximation in \autoref{eq9} for the BESS problem reduces to:

\begin{equation}
\label{eq41}
 \sum_{m\in\mathcal{M}}\sum_{t\in\mathcal{T}_0}P_{m,t}(c_{m,t}-d_{m,t})+\sum_{m\in\mathcal{M}}\sum_{t\in\mathcal{T}_\mathrm{obs}}\sum_{s\in\mathcal{S}}\pi_{s}P_{m,t,s}(c_{m,t,s}-d_{m,t,s}).
\end{equation}

\section{Results}
\label{sec5}

The formulation described in \autoref{eq9} is deployed on the systems outlined in \autoref{sec3} and \autoref{sec4} with $\alpha=0.95$. The intertemporal coupling constraints imposed by the storage processes make the scheduling problems non-myopic. In the IHS, the hydrogen inventory balance propagates the effect of production, storage, and consumption decisions into future periods. Meanwhile, the SOC dynamics impose the propagation for charge/discharge into the future BESS decisions. Consequently, decisions are path dependent: early market commitments determine the future state of storage assets and can limit the ability to make later actions and hedge against uncertainty. This is precisely the mechanism through which risk-averse sizing and dispatch can create value through the proposed stochastic framework. Our framework permits coordinated DA and ID participation rather than imposing strict market exclusivity. For BESS, exclusivity is only enforced between charging and discharging in the same market and period, while total activity is limited by inverter capacity. For IHS, participation is coupled through the process capacity and inventory constraints.

In both case studies, we restrict ourselves to the ``price-taker'' assumption, which can potentially result in suboptimal dispatch strategies \citep{gao2022multiscale} as well as under-sized units \citep{tsay2019optimal} for energy systems optimization at scale. As price-takers, we treat electricity prices as exogeneous stochastic processes. This assumption is appropriate for small- to medium-scale storage assets whose dispatch does not significantly influence the market, and it allows us to isolate the effect of uncertainty and downside risk management on design and operational decisions. For large-scale storage systems, however, dispatch may influence clearing prices and therefore feed back into future arbitrage opportunities. Under such price-maker or equilibrium-based formulations \mbox{\citep{tsay2023integrating,jalving2023beyond}}, aggressive dispatch strategies may have diminished profitability as these actions can erode the spread that they seek to exploit. At these large scales, market effects may in fact further motivate treating electricity prices as ``uncertain.'' Indeed, the stochastic formulation presented herein is exactly intended to enable dispatch adjustment upon changing market prices potentially induced by market participation. The present work does not model market-clearing conditions and coupling is therefore modelled through the shared operational decisions and physical constraints of the storage system; this modelling choice deliberately abstracts from system-level interaction to focus on the operational scheduling and risk-management problem from the operator perspective.

Although the case studies in this work focus on standalone systems integrated with electricity markets, some settings may require decision-making across multiple heterogenous assets (e.g., \mbox{\citet{banerjee2024sustainable}}). The CvaR-constrained formulation in \mbox{\autoref{eq9}} can be extended to these integrated multi-energy systems by embedding additional coupling constraints and imposing a suitable recourse structure for joint electricity, heat, and fuel systems. Importantly, the CVaR term would still quantify overall downside risk, while the process model would determine how the uncertainty propagates between energy domains and storage pathways. Quantitative behaviour of risk-aware sizing and scheduling may differ substantially in strongly coupled systems, particularly when integration creates alternative flexibility mechanisms beyond electricity arbitrage in a single storage medium.

We use traditional stochastic optimization metrics like expected value of perfect information (EVPI) and value of stochastic solution (VSS) \cite{birge1982value} to assess the impact of the proposed scheme on the respective systems. While these metrics do not explicitly account for the potential benefits of constraining tail risk, they can be used to quantify the sacrifice in expected cost induced by varying levels of risk aversion. To account for the benefit of constraining tail risk, we compute an adjusted $\mathrm{VSS}_{\mathrm{CVaR}}$ metric:

\begin{equation}
\label{eq42}
\mathrm{VSS}_{\mathrm{CVaR}}=\mathrm{VSS}+(\mathrm{CVaR_{SP_{\infty}}}-\mathrm{CVaR}_{\mathrm{SP}_{\epsilon}})+(E[\mathcal{L}_{\mathrm{SP}_{\infty}}]-E[\mathcal{L}_{\mathrm{SP}_{\epsilon}}]),
\end{equation} 

which penalizes losses induced on the expected cost by constraining CVaR, while also accounting for reductions in tail risk by computing the differences between a stochastic problem (SP) without an active CVaR constraint $\epsilon=\infty$ to those with a CVaR limit at $\epsilon$. To compute these metrics, we also distinguish between the expected solution of the expected value problem (EEV) (i.e., optimizing the first stage subject to the expected value of uncertainty and providing the first stage-solution to solve a second stage problem) and the wait-and-see problem (WS) (i.e., the expected solution having perfect information of the uncertainties). 

\begin{figure}[t]
\centering
\includegraphics[width=1\textwidth]{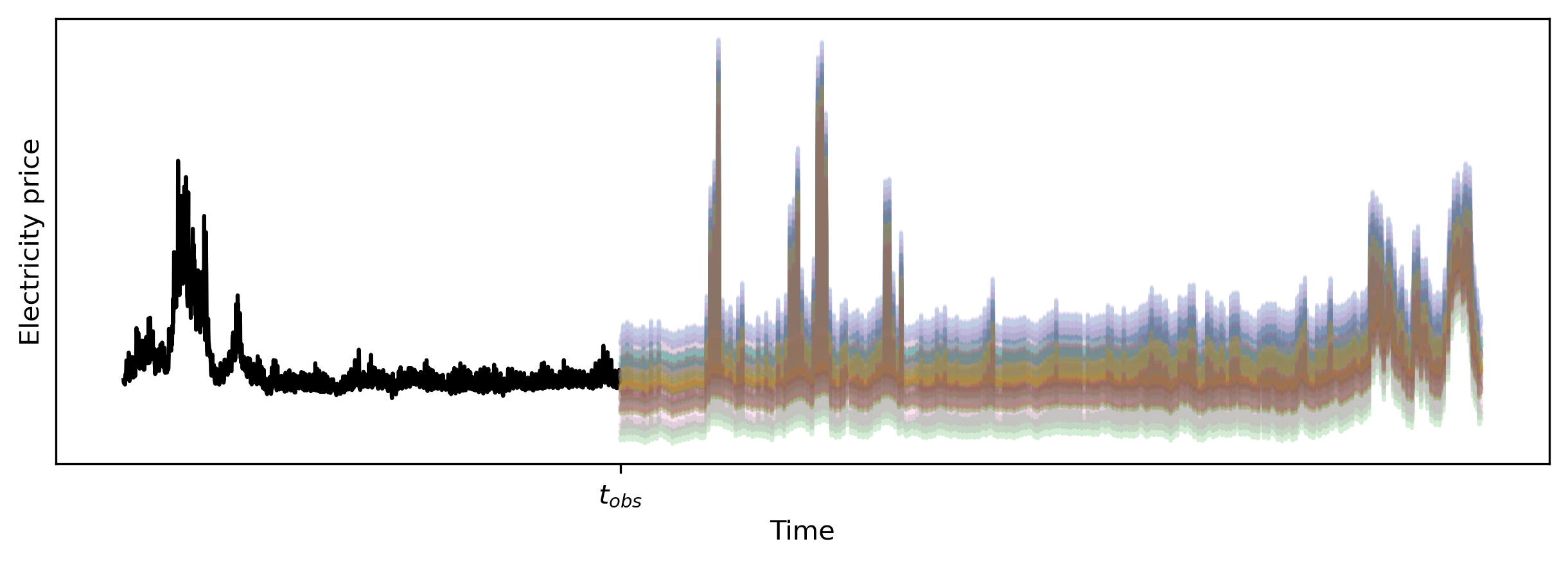}
\caption{Two-stage stochastic price structure - first-stage (initial trajectory) deterministic, second-stage (latter trajectory) uncertain.}\label{figextra}
\end{figure}

Both case studies are optimized over a yearly horizon, where an averaged price trajectory \mbox{$\mathbf{P_{t_{0},nom}}$} is used for the first stage and second-stage price trajectories are sampled from the distribution \mbox{$\xi\sim\mathcal{N}(\mathbf{P_{t_{obs},nom}},\sigma_\mathrm{obs}^2)$} as shown in \mbox{\autoref{figextra}}. We note that no explicit correlation between the markets was induced during the sampling procedure. However, the nominal price series used, which are described in the following sections, contain natural price correlations as they were retrieved for the same power authorities and time periods. This also highlights a benefit of using the general-purpose formulation in \mbox{\autoref{eq9}}, which only requires price series inputs rather than structural assumptions about uncertainty. As the uncertainty is limited to the electricity price, we calculate CVaR based only on operational cost (i.e., capital expenditure is deterministic); nonetheless, the formulation in \autoref{eq9} herein can admit uncertainty in both costs. Alternative scenario-generation methods (e.g., \citep{bounitsis2022data}) are also available in the literature. The scenario-based SAA used in this work introduces a trade-off between representation of uncertainty and computational tractability. As the number of scenarios, time periods, or recourse variables increases, the resulting optimization problem grows correspondingly. This scaling is particularly pronounced in the year-horizon, hourly resolution, and multi-market participation settings studied herein. For larger models and uncertainties with more complex probability distributions, strategies like Bender’s decomposition \mbox{\citep{garcia2023benders}}, scenario reduction \mbox{\citep{kim2023stochastic}}, temporal aggregation \mbox{\citep{teichgraeber2022time}}, and progressive hedging \mbox{\citep{qi2023portfolio}} may further improve tractability. Computational scaling was assessed systematically with respect to the number of scenarios used in the SAA (\mbox{\autoref{fig4}} and \mbox{\autoref{fig6}}). Expanding to more markets and/or multiple storage units could further increase the number of decision variables and model constraints, with corresponding scaling in computational effort.

All optimization problems and data analysis were performed on an Apple M3 Pro CPU. The optimization problems were solved using Gurobi 12.0.1 \mbox{\citep{gurobi}}; nominal computational times are reported with the scenario scaling analyses for each case study in the following sections.

\subsection{Yearly optimization of IHS}
\label{subsec14}

The process model $\mathbf{f}$, which imposes physical constraints for optimization of the IHS, is outlined in \autoref{sec3}. For this process, we aim to minimize capital and operating cost of the hydrogen production, storage, and conversion into electricity. Historical 2024 New England Independent Systems Operator (ISO-NE) hourly price data were retrieved from the \mbox{\citet{ISONEdata}} to be used as inputs for optimization. These represent the nominal electricity price trajectories \mbox{$\mathbf{P_{t_{0}, nom}}$} and  \mbox{$\mathbf{P_{t_{obs},nom}}$} denominated in \$US/kWh. We assume a value of $\sigma_\mathrm{obs}=\$20$/kWh as the price distribution standard deviation. This nominal and standard deviation parametrize the distribution of second-stage costs. The stochastic optimization of this system aims to optimally size the unit capacities at investment time, allowing for the operating decisions to be determined at a later time. The deterministic IHS design and scheduling problem has 164,011 variables and 180,405 constraints, scaling with the number of scenarios in the SAA.

\begin{figure}[t]
\centering
\includegraphics[width=1\textwidth]{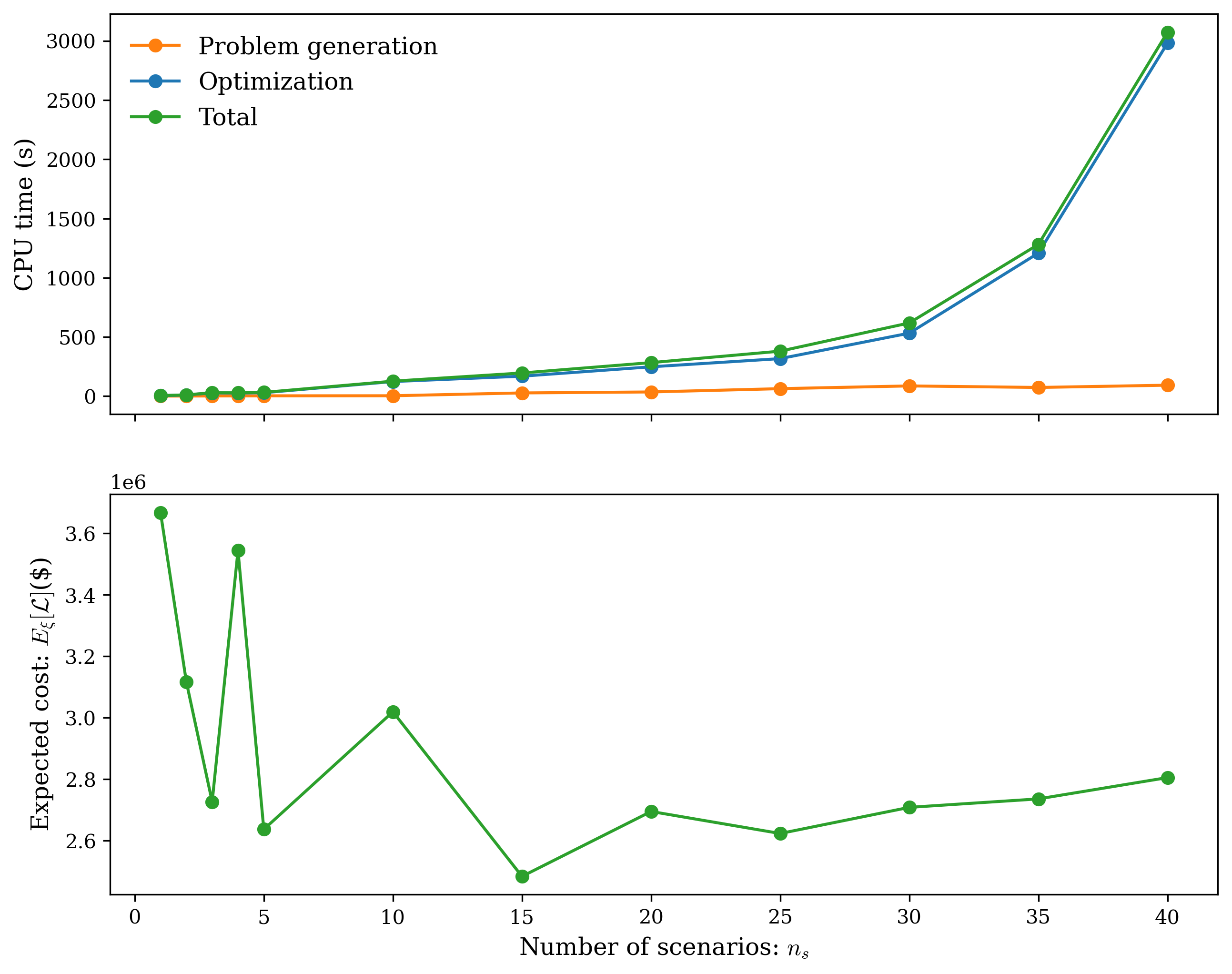}
\caption{Scaling of computational effort with discretization quality for IHS.}\label{fig4}
\end{figure}

We conducted a sensitivity analysis to understand the scaling between the number of scenarios, CPU time, and stability of solution afforded by the SAA. For this, we assume that true price trajectories can be observed after optimization of the design variables (i.e., $t_\mathrm{obs}=0$ hours). This observation setting corresponds to the case where the price trajectories are revealed upon completion of the plant build and at the beginning of operation of the IHS. Observation settings where $t_\mathrm{obs}>0$ correspond to cases where the price trajectories (i.e., whether prices are rising or falling) are not precisely known when plant operation begins. The former setting is used for sensitivity analysis as it produces the case with the most decision variables, and hence the most computational effort, as it requires no non-anticipativity constraints in the time domain. The results from this sensitivity analysis are displayed in \autoref{fig4}.

Based on the trade-offs between computation requirements and solution quality shown in \autoref{fig4}, we choose to formulate the SAA with $n_s = 35$ samples to balance computational load and approximation accuracy. In principle, many more scenarios can be used for a yearly optimization setting, as the solution does not need to be deployed online; however, this solution provides adequate solution stability to perform many optimization runs for comprehensive testing. 

We proceed to explore the trade-off between tail-risk, as measured by operating CVaR, and expected IHS cost by performing sensitivity analyses on the observation time $t_\mathrm{obs}$ and the CVaR bound $\epsilon$. \autoref{tab4} and \autoref{tab5} summarize these results. As shown in \autoref{tab4}, the solutions of all stochastic optimization problems, and that of the  EEV problem, do note include building a fuel cell. The former can be attributed to the fact that the stochastic program must hedge against scenarios in which there are low energy prices, hence performing hydrogen arbitrage would be economically unfavorable given the capital expenditure (i.e., the possibility of low prices induces sunk cost risk aversion). In contrast, the solution to the EEV problem does not include a fuel cell, as the nominal energy price timeseries does not contain high enough values to offset hydrogen production costs by reconverting to energy. These are contrasted to the case of perfect information (WS), where a small fuel cell is built on expectation as the potential gains from the known high price scenarios outweigh the costs from the low price scenarios. Note these analyses may change for more efficient fuel cell technologies, which could vary the parameters in \autoref{tab1}.

Interestingly, the solutions to all optimization formulations aside from the WS problem and the stochastic problems with $t_\mathrm{obs}=0$ hours include a compressor capacity equal to its lower bound. In the SPs, this occurs as there is large time-dependent uncertainty in the electricity prices as a single averaged price is not used for any portion of the trajectory when $t_\mathrm{obs}=0$ hours, hence the compressor must be over-sized to accommodate for the possibility of supplying the DRI furnace through the storage unit. As a fuel cell is built in the case of the WS case, more compression is required for supply to the cell.

\begin{table}[t]
  \centering
  \begin{tabular}{llllllll}
    \toprule
       $t_\mathrm{obs}$ & $\epsilon$ & $C_\mathrm{elec}$ & $C_\mathrm{stor}$ & $C_\mathrm{heat}$ & $C_\mathrm{comp}$ & $C_{fc}$ & $J^\mathrm{cap}$\\
    \midrule
         (hr) & (\$mn) &(GW) & (tn) & (MW) & (MW) & (MW) & (\$k) \\
    \midrule
    \multicolumn{2}{c}{WS} & $9.74$ & $0.226$ & $680$ & $142$ & $483$ & $414$\\
    \midrule
    \multicolumn{2}{c}{EEV} & $8.81$ & $0.171$ & $656$ & $100$ & $0$ & $353$\\
    \midrule
    \multirow{4}{*}{$0$} & $5.810^a$ & $9.08$ & $0.219$ & $723$ & $193$ & $0$ & $369$\\
    & $5.800$ & $9.92$ & $0.258$ & $734$ & $196$ & $0$ & $379$\\
    & $5.750$ & $10.35$ & $0.530$ & $823$ & $221$ & $0$ & $450$\\
    & $5.700$ & $11.92$ & $0.827$ & $947$ & $254$ & $0$ & $542$\\
    \midrule
     \multirow{4}{*}{$0^b$} & $7.050^a$ & $9.03$ & $0.223$ & $720$ & $100$ & $0$ & $367$\\
    & $7.000$ & $9.09$ & $0.237$ & $725$ & $100$ & $0$ & $371$\\
    & $6.950$ & $10.22$ & $0.511$ & $814$ & $100$ & $0$ & $443$\\
    & $6.900$ & $11.76$ & $0.822$ & $937$ & $100$ & $0$ & $534$\\
    \midrule
    \multirow{4}{*}{$2200$} & $5.815^a$ & $8.93$ & $0.196$ & $708$ & $100$ & $0$ & $360$\\
    & $5.800$ & $9.19$ & $0.266$ & $727$ & $100$ & $0$ & $378$\\
    & $5.750$ & $10.29$ & $0.544$ & $814$ & $100$ & $0$ & $449$\\
    & $5.700$ & $11.86$ & $0.843$ & $938$ & $100$ & $0$ & $540$\\
    \midrule
    \multirow{4}{*}{$4400$} & $5.820^a$ & $8.84$ & $0.178$ & $682$ & $100$ & $0$ & $355$\\
    & $5.800$ & $9.18$ & $0.270$ & $702$ & $100$ & $0$ & $378$\\
    & $5.750$ & $10.27$ & $0.552$ & $773$ & $100$ & $0$ & $449$\\
    & $5.700$ & $11.79$ & $0.846$ & $876$ & $100$ & $0$ & $540$\\ 
    \bottomrule
  \end{tabular}
  \caption{Expected design and capital decisions for varying IHS stochastic program hyperparameters. $^a$ denotes that $\epsilon$ was not bounded (i.e., the solution is risk-neutral), otherwise the solution occurs at the risk bound (i.e., CVaR$=\epsilon$).}\label{tab4}
\end{table}

Comparing across the stochastic problems, there is large upward monotonicity in the electrolyzer, storage, and heat capacities with increasing risk aversion (decreasing $\epsilon$). For instance, tightening the CVaR bound from $\epsilon = \infty$ to $\epsilon = \$ 5.7$mn at $t_0=0$ increased electrolyzer capacity by $\approx 31 \% $, storage capacity by $\approx 278 \% $, and heater capacity by $\approx 31 \% $. This trend suggests that building overcapacity is the more risk-averse design policy for the IHS system. As CVaR is more constrained, the solutions to the stochastic program are more conservative, including increasing capacities to accommodate for extreme scenarios that could lead to high losses. Crucially, this overcapacity allows for larger hydrogen stores and heating supply to ensure that the DRI demand is met at all times. These overcapacities are all further reflected in increasing capital cost with increasing CVaR bound.

The observation time has a smaller, but opposite, effect to the risk aversion level. Increasingly delayed observations of the true price trajectories cause lower capacities to be built; this occurs as an averaged electricity price trajectory is used for the first stage, which the stochastic program takes as the ``true'' price. The use of a sample-averaged trajectory over long periods of time forces the solution of the stochastic program to include capacity that optimizes the sample-averaged objective, rather than to hedge against potential uncertainties by building overcapacity. The suboptimalities induced by this approximation are reflected in the operating costs outlined in \autoref{tab5}.

\begin{table}[!ht]
  \centering
  \begin{tabular}{lllllll}
    \toprule
       $t_\mathrm{obs}$ & $\mathrm{CVaR}$ & $E_{\xi}[\mathcal{L}]$ & $\frac{E_{\xi}[J^{op}_{ID}]}{E_{\xi}[J^{op}_{DA}]}$ & EVPI & VSS & $\mathrm{VSS_{CVaR}}$\\
    \midrule
         (hr) & (\$mn) & (\$mn) & $-$ & (\$k) & (\$k) & (\$k) \\
    \midrule
     & WS & $2.716$ & $1.278$ & $-$ & $-$ & $-$\\
    \midrule
    \multirow{5}{*}{$0$} & EEV  & $2.799$ & $1.284$ & $83$ & $-$ & $-$\\
    & $5.810^a$ & $2.735$ & $1.299$ & $19$ & $63$ & $63$\\
    & $5.800$ & $2.736$ & $1.308$ & $20$ & $63$ & $72$\\
    & $5.750$ & $2.748$ & $1.363$ & $32$ & $50$ & $97$\\
    & $5.700$ & $2.776$ & $1.412$ & $60$ & $23$ & $92$\\
    \midrule
    \multirow{5}{*}{$0^b$} & EEV  & $2.899$ & $1.289$ & $183$ & $-$ & $-$\\
    & $7.050^a$ & $2.870$ & $1.302$ & $71$ & $29$ & $29$\\
    & $7.000$ & $2.870$ & $1.306$ & $71$ & $29$ & $79$\\
    & $6.950$ & $2.883$ & $1.364$ & $167$ & $16$ & $103$\\
    & $6.900$ & $2.913$ & $1.419$ & $197$ & $-14$ & $93$\\
    \midrule
    \multirow{5}{*}{$2200$} & EEV  & $2.799$ & $1.283$ & $83$ & $-$ & $-$\\
    & $5.815^a$ & $2.736$ & $1.290$ & $21$ & $62$ & $62$\\
    & $5.800$ & $2.738$ & $1.308$ & $22$ & $61$ & $74$\\
    & $5.750$ & $2.753$ & $1.362$ & $37$ & $46$ & $95$\\
    & $5.700$ & $2.785$ & $1.418$ & $69$ & $14$ & $81$\\
    \midrule
    \multirow{5}{*}{$4400$} & EEV  & $2.800$ & $1.283$ & $83$ & $-$ & $-$\\
    & $5.820^a$ & $2.738$ & $1.289$ & $21$ & $62$ & $62$\\
    & $5.800$ & $2.740$ & $1.308$ & $24$ & $59$ & $77$\\
    & $5.750$ & $2.757$ & $1.364$ & $41$ & $42$ & $92$\\
    & $5.700$ & $2.792$ & $1.419$ & $76$ & $7$ & $73$\\
    \bottomrule
  \end{tabular}
  \caption{Operating costs and stochastic summary metrics for varying IHS stochastic program hyperparameters. $^a$ denotes that $\epsilon$ was not bounded (i.e., the solution is risk-neutral), otherwise the solution occurs at the risk bound (i.e., CVaR$=\epsilon$).}\label{tab5}
\end{table}

As shown in \autoref{tab5}, increasingly delayed observations correspond to larger expected costs. Intuitively, the earlier one can access true price information, the better decisions one can take. 
This results in lower VSS (i.e., lower benefits in optimizing stochastically) but higher EVPI (i.e., more solution suboptimality with respect to the wait-and-see problem). 

\begin{figure}[t]
\centering
\includegraphics[width=1\textwidth]{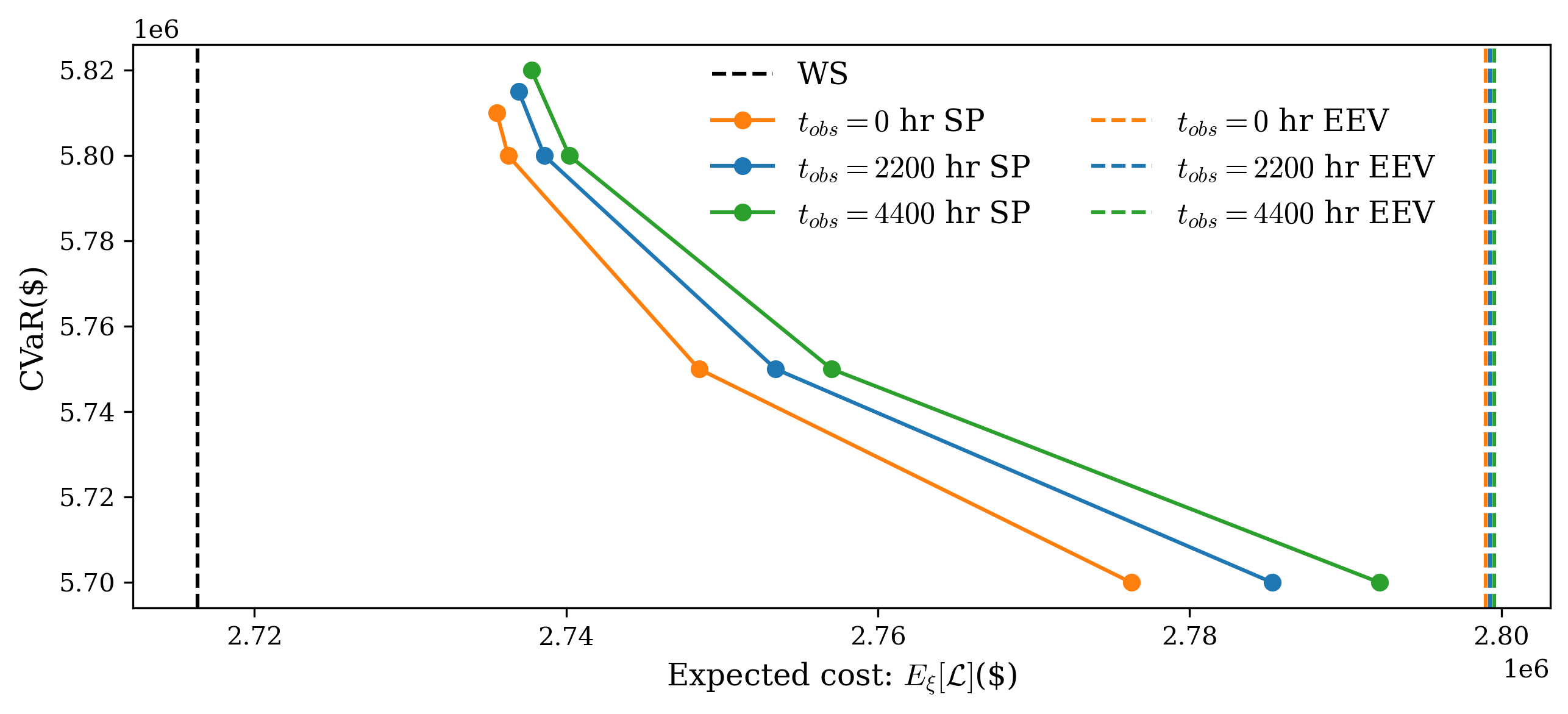}
\caption{Trade-off between expected cost and risk aversion for IHS. }\label{fig5}
\end{figure}

We found the ratio of ID to DA costs to increase with risk aversion (lower $\epsilon$). This may be attributed to the over-sized capacities shown in \autoref{tab4}, which enable more purchasing and subsequent storage from the volatile ID market. Particularly in periods of low ID prices, these large capacities enable increased storage to buffer future DRI energy requirements; this corresponds to a more risk-averse strategy as the storage system is effectively hedging the effects of potentially high future prices to satisfy system power demand. The trade-off between risk aversion and expected cost is also shown in \autoref{fig5} for varying $t_\mathrm{obs}$.

As shown in \autoref{fig5}, there is a risk-reward trade-off between CVaR and expected IHS cost. Tighter CVaR bounds increase the expected cost of the system while limiting extreme losses. In general, the range of expected cost sacrifice is $ \approx \$60$k while the potential reductions in CVaR have a range of $\approx\$115$k; hence the potential benefits outweigh the cost. The lowest expected cost occurs when CVaR is not constrained (i.e., $\epsilon=\infty$); this results in $\mathrm{CVaR} \approx \$5.85$mn beyond which expected cost cannot be reduced to the level exhibited by the wait-and-see (perfect information) solution; this gap represents the EVPI. The largest expected costs occur when the CVaR bound is tightest (i.e., $\epsilon=\$5.7$mn), beyond which CVaR cannot be reduced. Nevertheless, these solutions with the tightest risk bounds outperform the EEV solutions; this gap represents the VSS. The benefits of modeling tail risk are evident in the $\mathrm{VSS_{CVaR}}$ as $\mathrm{VSS_{CVaR} \geq VSS}$ for all risk-constrained scenarios; this means that the tail risk avoided is greater than the respective increase in expected cost incurred. That is, there is always an outsized benefit in the trade-off betweeen expected value and CVaR for the IHS case study.

The IHS case study was also tested under second-stage price trajectories sampled from a Student’s t-distribution with $\nu=5$ and a scaling factor of $\sigma=15.55$. These parameters provide the t-distribution with equivalent variance to the normal distribution used throughout this work, while effectively providing an ablation study with a heavy-tailed distribution. These results can be found in \mbox{\autoref{tab4} and \autoref{tab5}}, where the heavy tails result in CVaR values increased by $\approx \$1.1$mn with respect to the normal distribution. Despite this more extreme setting, the risk-constrained optimization approach exhibits equivalent sensitivity to the CVaR bound for both distributions, as shown by the $\mathrm{VSS_{CVaR}}$ values. Indeed, the optimizer can only use the capacities to hedge against risk in the IHS plant, and this leads to similar risk abatement behaviour regardless of the uncertainty distribution. Although the risk-constrained formulation is found beneficial regardless of the weight of the tail, we also note that the EVPI is higher for the t-distribution cases; this suggests the SAA approximation is worse when the uncertainty distribution is less straightforward.

\subsection{Yearly optimization of BESS}
\label{subsec15}

The process model $\mathbf{f}$, which imposes physical constraints for optimization of the BESS, is outlined in \autoref{sec4}. The deterministic BESS dispatch problem has 245,280 variables and 3592,58 constraints, scaling with the number of scenarios in the SAA. For this process, we aim to maximize the operating profit of the energy arbitrage enabled by battery storage. Price data were generated based on industry-acquired predictions to be used for optimization; these represent the nominal electricity price trajectories $\mathbf{P_{t_{0}, nom}}$ and $\mathbf{P_{t_{obs},nom}}$ denominated in \euro$/kWh$. This yearly optimization setting corresponds to various potential power purchase agreements (PPAs) according to the observation time $t_\mathrm{obs}$. Specifically, $t_\mathrm{obs}=0$ corresponds to a fixed PPA where the electricity prices are fixed at agreement time; conversely, the flexible PPA is modelled as a contract in which electricity prices are not fully fixed at agreement time but are revealed after a delay $t_\mathrm{obs}>0$. This captures practical PPA structures in which a buyer commits to volume and availability, while the effective price is indexed to market outcomes or a moving reference. Within the two-stage stochastic framework, recourse decisions allow the BESS operator to adapt dispatch once prices are revealed. The proposed risk-constrained formulation directly supports the design and operation of flexible PPAs by enabling operators to quantify downside risk exposure associated with uncertain contract pricing, select PPA structures (e.g., shorter vs longer observation times) that satisfy explicit tail-risk limits via CVaR constraints, and determine operational strategies (DA vs ID participation) that hedge against unfavorable realizations of flexible pricing terms. In this sense, the model provides a decision-support tool for negotiating PPA flexibility clauses and for selecting dispatch strategies that remain profitable under worst-case pricing outcomes.

\begin{figure}[t]
\centering
\includegraphics[width=1\textwidth]{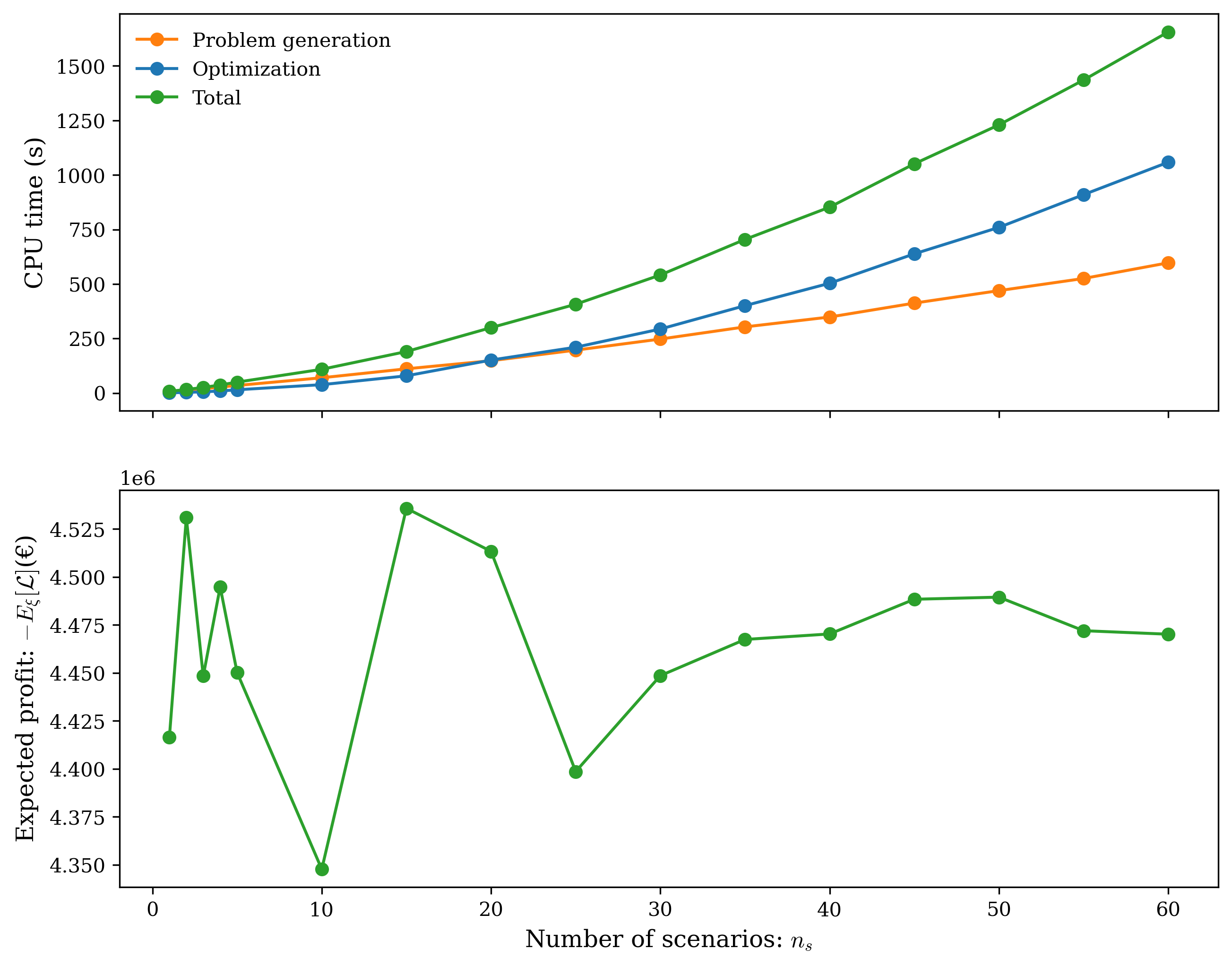}
\caption{Scaling of computational effort with discretization quality for BESS.}\label{fig6}
\end{figure}

Similarly to the IHS, we first conduct a sensitivity analysis to understand the scaling between the number of scenarios, CPU time, and stability of solution afforded by the SAA. $\sigma_\mathrm{obs}=\;$\euro$\;30$/kWh is assumed. The results from this sensitivity analysis are shown in \autoref{fig6}. Based on the scaling of computation requirements and solution quality, we again choose to formulate our SAA using $n_s = 35$ samples to balance computational load and approximation accuracy. Since this system is subject to the same uncertainties (i.e., energy prices) and solves the same energy allocation problem (albeit in a different system), the scaling of computational effort follows a similar trend as in the IHS (\autoref{fig4}).

We examine the trade-off between the tail risk and expected arbitrage profit under various choices for the values of optimization hyperparameters, as summarized in \autoref{tab6}. We note that the expected costs are denoted with a negative sign (i.e., $-E_{\xi}[\mathcal{L}]$) as their values correspond to profits, following the minimization convention in \autoref{eq9}, which is common in the literature. For the BESS case study, we again observe the expected profit to be more sensitive to the risk bound than to observation time. Tightening risk bounds (i.e., constraining potential losses) and delaying observation times again both result in decreasing expected profits. The tail risks (CVaR) and risk bounds ($\epsilon$) all fall in the non-negative range, as the tail risk constitutes overall loss or break-even scenarios, respectively. With a tail risk range of $\approx\;$\euro$\;1.25$mn, using a risk-constrained approach can effectively eliminate the possibility of losses by sacrificing $\approx\;$\euro$\;506$k in the lowest profit scenario (i.e., a late observation time). 

\begin{table}[!ht]
  \centering
  \begin{tabular}{lllllll}
    \toprule
       $t_\mathrm{obs}$ & $\mathrm{CVaR}$ & $-E_{\xi}[\mathcal{L}]$ & $\frac{E_{\xi}[J^{op}_{ID}]}{E_{\xi}[J^{op}_{DA}]}$ & EVPI & VSS & $\mathrm{VSS_{CVaR}}$\\
    \midrule
         (hr) & (\euro mn) & (\euro mn) & $-$ & (\euro k) & (\euro k) & (\euro k) \\
    \midrule
     & WS & $4.467$ & $293$ & $-$ & $-$ & $-$\\
    \midrule
    \multirow{7}{*}{$0$} & EEV  & $4.467$ & $299$ & $0.5$ & $-$ & $-$\\
    & $1.25^a$ & $4.467$ & $293$ & $0$ & $0.5$ & $0.5$\\
    & $1.00$ & $4.465$ & $469$ & $2.4$ & $-1.9$ & $245.8$\\
    & $0.75$ & $4.455$ & $-419$ & $12.4$ & $-11.9$ & $475.7$\\
    & $0.50$ & $4.424$ & $-81$ & $43.8$ & $-43.3$ & $662.9$\\
    & $0.25$ & $4.343$ & $-36$ & $124.5$ & $-124.0$ & $751.5$\\
    & $0.00$ & $4.112$ & $-16$ & $355.2$ & $-354.7$ & $540.0$\\
    \midrule
    \multirow{7}{*}{$2200$} & EEV  & $4.464$ & $296$ & $3.7$ & $-$ & $-$\\
    & $1.25^a$ & $4.466$ & $284$ & $1.8$ & $2.0$ & $2.0$\\
    & $1.00$ & $4.461$ & $457$ & $5.6$ & $-1.9$ & $244.3$\\
    & $0.75$ & $4.445$ & $-189$ & $21.9$ & $-18.1$ & $461.8$\\
    & $0.50$ & $4.402$ & $-51$ & $65.0$ & $-61.2$ & $625.4$\\
    & $0.25$ & $4.302$ & $-24$ & $165.1$ & $-161.3$ & $675.3$\\
    & $0.00$ & $4.046$ & $-14$ & $421.3$ & $-417.5$ & $412.9$\\
    \midrule
    \multirow{7}{*}{$4400$} & EEV  & $4.457$ & $411$ & $10.8$ & $-$ & $-$\\
    & $1.25^a$ & $4.463$ & $440$ & $4.6$ & $6.2$ & $6.2$\\
    & $1.00$ & $4.456$ & $898$ & $11.0$ & $-0.3$ & $243.3$\\
    & $0.75$ & $4.432$ & $-113$ & $35.4$ & $-24.6$ & $444.6$\\
    & $0.50$ & $4.375$ & $-40$ & $92.0$ & $-81.2$ & $581.3$\\
    & $0.25$ & $4.257$ & $-20$ & $210.2$ & $-199.4$ & $595.9$\\
    & $0.00$ & $3.961$ & $-13$ & $506.1$ & $-495.3$ & $253.1$\\
    \bottomrule
  \end{tabular}
  \caption{Operating costs and stochastic summary metrics for varying BESS stochastic program hyperparameters. $^a$ denotes that $\epsilon$ was not bounded (i.e., the solution is risk-neutral), otherwise the solution occurs at the risk bound (i.e., CVaR$=\epsilon$).}\label{tab6}
\end{table}

We note that the market participation ratio is negative in some settings, corresponding to when losses are incurred in the DA market and all profits are made from participation in the ID market. The majority of the expected profit comes from the ID  market, because its large volatility produces the best buying and selling opportunities for arbitrage. However, the participation in the DA market increases with tightening risk bound as more overall energy is bought to increase the inventory over time. Intuitively, balancing participation in the more aggressive ID market and the more conservative DA market enables controlling the risk-reward tradeoff.

From a stochastic optimization perspective, the expected value of perfect information (EVPI) is very sensitive to the risk bound as the risk tolerance induces large gaps in profit with respect to the wait-and-see (WS) problem. Furthermore, the value of stochastic solution (VSS) in this case can be negative as the expected solution of the expected value problem (EEV) outperforms the risk-constrained stochastic problem (SP); this is owed to the expected profit reduction induced when constraining CVaR. However, when considering the potential for risk avoidance, the benefits for constraining risk as reflected in $\mathrm{VSS_{CVaR}}$ are again significant, with maximum benefits of $\approx\;$\euro$\;750$k. Similarly to in the IHS case study, the potential losses incurred by limiting risk are outweighed by the benefit according to this metric.

\begin{figure}[t]
\centering
\includegraphics[width=1\textwidth]{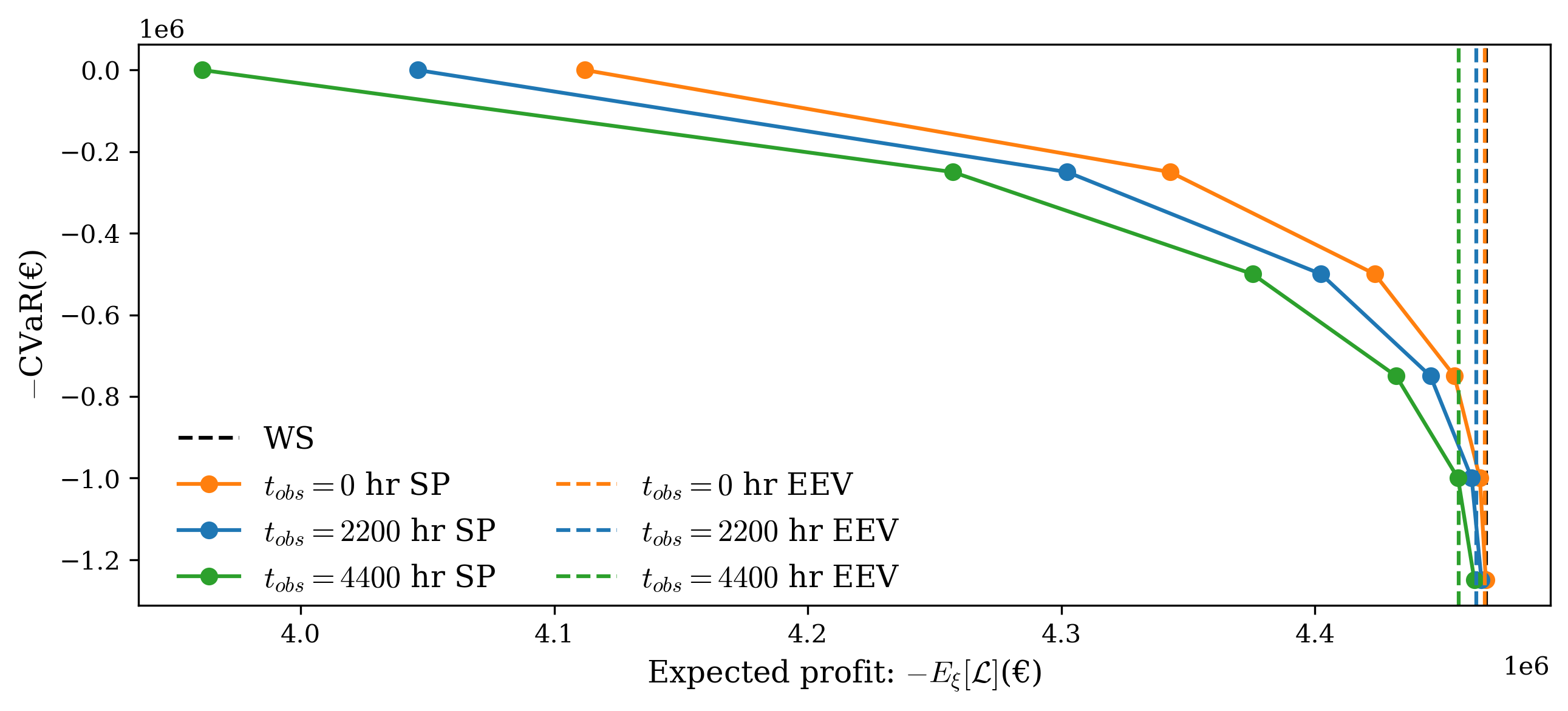}
\caption{Trade-off between expected cost and risk aversion for BESS. }\label{fig7}
\end{figure}

The risk-reward trade-off for the BESS is shown in \autoref{fig7}. All EEV, unbounded SP, and WS solutions have similar expected profits; this is reflective of the relatively low VSS and EVPI. Conversely, the large range of CVaR can also be observed, which reflects the high risk profile of the BESS arbitrage problem. Irrespective of the observation time, the CVaR can always be reduced to zero, which constitutes a break-even scenario. The ability constrain tail risk from a loss regime to a break-even regime is especially powerful for risk-averse operators, which elucidates the benefits afforded by CVaR-constrained optimization.

\begin{figure}[t]
\centering
\includegraphics[width=1\textwidth]{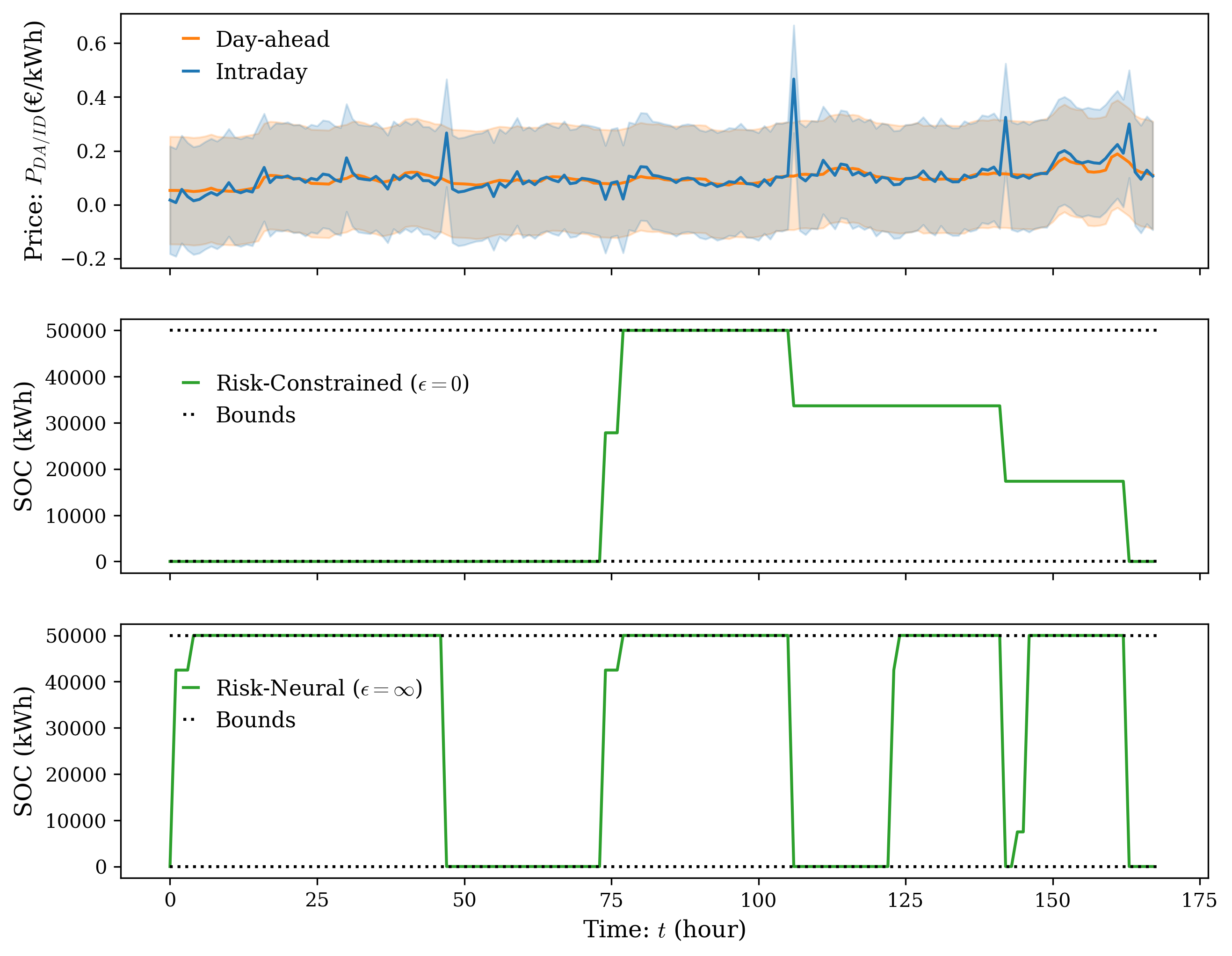}
\caption{Representative week for BESS schedule in with risk-constrained (middle) and risk-neutral (bottom) optimizers. Corresponding nominal prices and their associated uncertainty (top).}\label{ExtraFig2}
\end{figure}

Tightening the CVaR bound from $\epsilon = \infty$ to $\epsilon =$ \euro 0  at $t_0=0$ reduced expected profit by $\approx 8 \%$ while altogether eliminating tail loss (i.e., $100 \%$ reduction). This is achieved through less aggressive trading when risk is constrained as illustrated in Figure 9, which compares BESS scheduling under DA and ID price uncertainty and illustrates the impact of CVaR-based risk constraints on operational decisions. The top panel shows the expected prices together with their respective uncertainty envelopes, representing the range of plausible price realizations across scenarios. The ID price exhibits substantially wider dispersion and more pronounced tail events than the DA price, highlighting the presence of significant downside risk when arbitraging based on intraday forecasts. Under the most tight CVaR-constraint formulation with $\epsilon = 0$ (\autoref{ExtraFig2}, middle panel), the BESS follows a conservative charging policy. The SOC remains at zero for extended periods and increases sharply only when the expected price spread is sufficiently robust across the uncertainty set. This behaviour reflects the CVaR constraint, which limits the expected loss in the worst $\alpha$-percentile of scenarios and therefore discourages actions that are profitable on average but expose the operator to adverse tail outcomes. Consequently, the battery performs more gradual and monotonic charge-discharge cycles that avoids frequent reactions to transient ID price spikes that lie within the uncertainty envelope. This policy prioritizes robustness to unfavorable price realizations and implicitly reduces cycling-related degradation.

In contrast, the risk-neutral policy $\epsilon = \infty$ (\autoref{ExtraFig2}, bottom panel) ignores tail-risk exposure and optimizes based on expected profit. The resulting SOC trajectory exhibits rapid transitions between empty and full capacity, indicating aggressive exploitation of short-term price spreads and ID spikes, including those near the edges of the uncertainty evelope. While this strategy achieves higher expected arbitrage revenue, it is highly sensitive to forecast errors and tail events in the ID price distribution. In adverse realizations, the same actions may lead to substantial losses, which are explicitly controlled in the CVaR-constrained formulation but left unmitigated in the risk neutral case.

Overall, \autoref{ExtraFig2} exemplifies how incorporating CVaR constraints yields materially different operational behaviour even under identical expected price signals. The CVaR-constrained policy sacrifices short-term arbitrage opportunities in periods of high price uncertainty to hedge against tail losses, whereas the risk-neutral policy exploits volatility more aggressively at the cost of increased financial and operational risk. These results highlight our framework as an effective tuning mechanism to balance between expected profit and robustness for BESS participation in volatile ID markets.

\subsection{Rolling horizon optimization of BESS}
\label{subsec16}

\begin{figure}[t]
\centering
\includegraphics[width=1\textwidth]{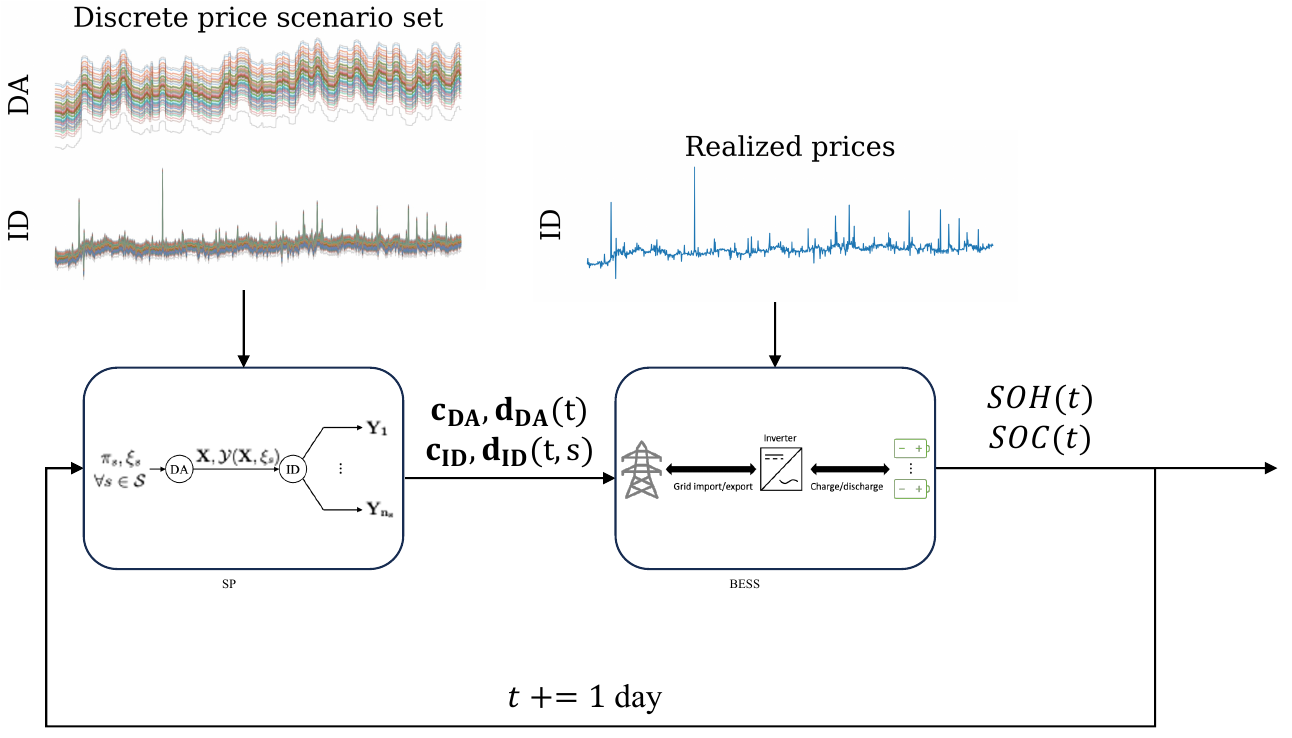}
\caption{Stochastic rolling horizon operation of BESS schedule. }\label{fig8}
\end{figure}

As the BESS system does not require any capacity decisions, a risk-constrained formulation can also be deployed in a rolling horizon manner for online operation. As depicted in \autoref{fig8}, a rolling horizon approach can continually re-optimize the charge/discharge policy of the system at a fixed time interval. By performing this re-optimization, updated price projections that comprise the scenario set can be incorporated into the stochastic optimization problem. Previous works also support that rolling-horizon optimization itself helps mitigate problem uncertainty~\citep{lejarza2022feedback,wang2022multi, mcallister2022inherent, risbeck2019economic}. A rolling horizon stochastic optimization approach has been applied to energy storage systems previously \citep{kumar2019benchmarking} and can leverage market timings where DA prices are set through auction in advance of the energy actually being deployed while the intraday prices are determined in a spot market with a central tendency around the previously-set DA price. 

This online setting corresponds to a live energy trading case in which DA comitments are made before DA prices are determined and ID prices follow a noisy distribution around the DA trajectories. Following the rolling-horizon approach, we forego the temporal partitioning of here-and-now and wait-and-see variables presented in \autoref{sec4} in order to exploit the market structure. The risk-constrained two-stage problem for this approach has the first-stage decisions $\mathbf{X} = \begin{bmatrix} \mathbf{c_{DA}} & \mathbf{d_{DA}} \end{bmatrix}^{\top}$,  where $\mathbf{c_{DA}} \in\mathcal{C}_{DA}\subset\mathbb{R}^{|\mathcal{T}|}$ are the charge and $\mathbf{d_{DA}}\in\mathcal{D}_{DA}\subset\mathbb{R}^{| \mathcal{T}| }$ are the discharge dispatch from the grid. This formulation determines a single dispatch for the DA actions that is optimal for the whole scenario set for as the prices have yet to be realized by auction. The second stage decisions are the charge/discharge dispatch decision in the intraday market $\mathbf{Y}=\begin{bmatrix} \mathbf{c_{ID}} & \mathbf{d_{ID}} \end{bmatrix}^{\top}$ where $\mathbf{c_{ID}}\in\mathcal{C}_{ID}(\mathbf{c_{DA}}, \mathbf{d_{DA}})\subset\mathbb{R}^{| \mathcal{T}| }$ and $\mathbf{d_{ID}}\in\mathcal{D}_{ID}(\mathbf{c_{DA}}, \mathbf{d_{DA}})\subset\mathbb{R}^{| \mathcal{T}|}$. Noticing that, once the DA prices are revealed, the ID prices follow their general trend (albeit noisily), we assume that the decisions corresponding to one of the scenarios in the scenario set are implemented as recourse. Accordingly, the first-stage price vector is $\mathbf{c}=\mathbf{P_{DA}}\in\mathbb{R}^{| \mathcal{T}|}$ and second-stage uncertain price vector is the multivariate distribution $\mathbf{\xi}=\mathbf{P_{ID}}\in\Psi\subset\mathbb{R}^{|\mathcal{T}|}$. As implied by the dimensionality of the variables, the time horizon $\mathcal{T}$ is used, which corresponds to a year-long window. The charge and discharge decisions for each day are further constrained according to the previous day's decisions by continually updating the $SOC$ and $SOH$ of the model. When discretized in time, the overall sample-average approximation in \autoref{eq9} for the BESS rolling horizon problem reduces to:

\begin{equation}
\label{eq43}
\sum_{t\in\mathcal{T}}P_{DA,t}(c_{DA,t}-d_{DA,t})+\sum_{t\in\mathcal{T}}\sum_{s\in\mathcal{S}}\pi_{s}P_{ID,t,s}(c_{ID,t,s}-d_{ID,t,s}).
\end{equation}

We perform the rolling horizon optimization of a single week of BESS operation. As displayed in \autoref{tab7}, we report cumulative metrics to account for each day optimized, and we impose risk constraints ($\epsilon$) that correspond to annualized values since a year-long window is used. In contrast to the year-long formulation without feedback, not all of the solutions to weekly optimization problems solved in the weekly operating period can limit the risk to a break-even ($\epsilon  =0$) setting; hence, this row is omitted from \autoref{tab7}. For the week considered, the optimal solutions do not involve participation in the DA market at all, and we therefore omit the market participation ratio. This corroborates the previous results in \autoref{subsec15} where the optimal results involve very little participation in the DA market. In the current setting, where participation in the DA market is optimized as a here-and-now decision, the stochastic program chooses the best policy is to wholly adapt its dispatch to a given charge level that corresponds to an ID trajectory centered around the DA auction prices.

\begin{table}[!ht]
  \centering
  \begin{tabular}{llllll}
    \toprule
       $\epsilon$ & $\sum\mathrm{CVaR}$ & $-\sum E_{\xi}[\mathcal{L}]$ & $\sum$ EVPI & $\sum$ VSS & $\sum\mathrm{VSS_{CVaR}}$\\
    \midrule
         (\euro mn) & (\euro k) & (\euro k) & (\euro k) & (\euro k) & (\euro k) \\
    \midrule
    \multicolumn{2}{c}{WS} & $0.1$ & $-$ & $-$ & $-$\\
    \midrule
    \multicolumn{2}{c}{EEV}  & $-35.2$ & $35.3$ & $-$ & $-$\\
    \midrule
    {$\infty$} & {$32$} & $-2.6$ & $2.7$ & $32.6$ & $32.6$\\
    {$1$} & {$28$} & $-2.5$ & $2.6$ & $32.7$ & $32.8$\\
    {$0.75$} & {$21$} & $-2.5$ & $2.6$ & $32.7$ & $32.8$\\
    {$0.5$} & {$19$} & $-1.9$ & $2.0$ & $33.2$ & $33.9$\\
    {$0.25$} & {$14$} & $-0.9$ & $1.0$ & $34.2$ & $35.9$\\
    \bottomrule
  \end{tabular}
  \caption{Operating costs and stochastic summary metrics for varying BESS rolling horizon stochastic program hyperparameters. Sums ($\sum$) denote cumulative quantities for a week of operation.}\label{tab7}
\end{table}

\begin{figure}[t]
\centering
\includegraphics[width=1\textwidth]{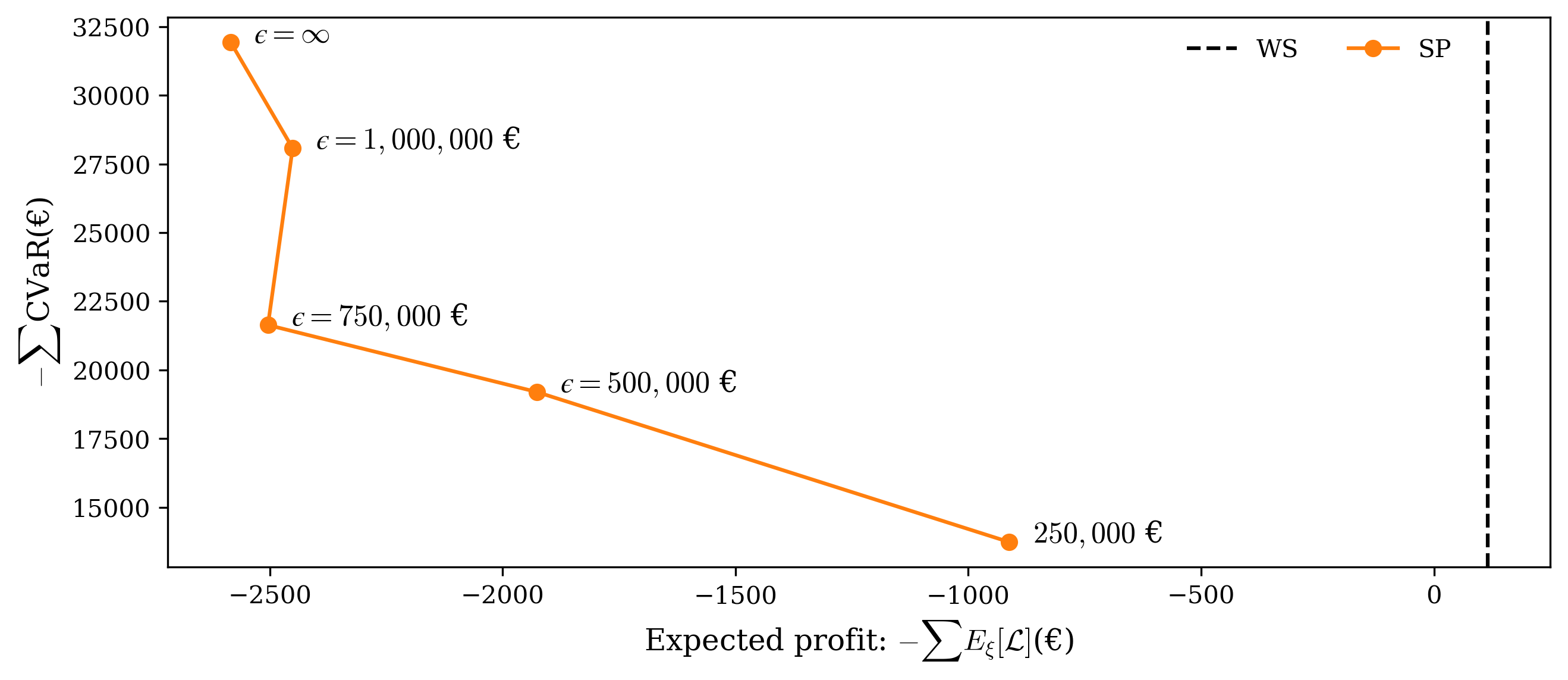}
\caption{Trade-off between expected cost and risk aversion for BESS rolling horizon operation. }\label{fig9}
\end{figure}

As illustrated in \autoref{fig9} and summarized in \autoref{tab7},  the risk-constrained formulation outperforms the EEV solution, while being less suboptimal to the WS solution (i.e., a high VSS and a low EVPI). In general, most of the VSS comes from the poor performance of the EEV problem, where the $\mathrm{VSS_{CVaR}}$ is nearly equivalent. However, with longer operating periods, we expect for the performance of the EEV to average out and for the SP to provide more benefit owing to risk avoidance as in \autoref{subsec15}. Counterintuitively, the expected cost in this case is lower with tighter risk bounds. This is potentially a result of the week being optimized, where a more aggressive charging policy in fact constitutes the most risk-averse strategy. \autoref{fig9} also shows some non-monotonicity in the expected profit as the risk bound approaches the non-constrained region. This occurs as limiting the CVaR does not guarantee overall monotonicity in the time domain, where several risk-averse problems are performed iteratively. In other words, the most conservative closed-loop (rolling horizon) solution does not necessarily correspond to the most conservative open-loop solution. Nevertheless, this operational scheme exhibits closer performance to the perfect information (WS) solution despite the presence of uncertainty while also providing shorfall-limiting potential.

\section{Conclusion}
\label{sec6}

A risk-constrained stochastic approach was proposed for energy systems optimization under electricity market uncertainty, with applications in the design and scheduling of a broad class of energy storage media. The proposed framework exposes how storage-specific state dynamics amplify tail-risk in multi-market participation, and how CVaR constraints reshape optimal DA–ID allocation strategies. Constraining CVaR is especially beneficial with respect to risk-neutral optimization in settings where uncertainty is high, penalties or reliability constraints are severe, and recourse is limited (e.g., two-stage optimisation). The risk-constrained approach allows for explicit specification of the operator's tolerance for risk rather than using a heuristic objective weight. This approach was tested in IHS and BESS case studies; the former includes capacity and dispatch decisions while the latter only involves dispatch but was also applied in a rolling horizon manner. Both systems were subject to DA and ID electricity markets whereby the optimizer determined the charge and discharge at each point in time to each market. This variety of settings across our case studies elucidates the flexibility of our framework to address uncertainty in scheduling for integrated design and operation, fixed or flexible PPA, and live trading uses.

Both systems exhibit a trade-off between the expected cost/profit and the tail risk. In the IHS cost minimization setting, higher expected costs resulted in an outsized constraining of the CVaRs. Further, a weaker correlation with the observation time of the true price signal was observed. The first-stage decisions are manifested in the IHS system through the unit sizes whereby large capacities were built for tighter risk bounds. This comes at the expense of capital cost but allows for larger inventories to be held at any given time; hence more participation in the ID market to take advantage of large price spreads. The relationship between expected operational profit and tail risk was evident to a larger extent on the BESS system, where profit was sacrificed to limit potential tail losses. Indeed, using an optimal risk-constrained schedule on the BESS system can shift the CVaR for a loss to a break-even regime. As in the IHS, the BESS arbitrage was observed to make an increasing amount of profit from the ID market with increasing risk aversion; this is reflected in order-of-magnitude larger ID profits and, in some cases, net losses to the DA market. As larger price spreads occur owed to ID volatility, this large proportion of ID market profit is appropriate in the arbitrage context. In a rolling horizon implementation, the stochastic optimization of the BESS system resulted in significant improvements over the expected value problem with only small deterioration with respect to the perfect information problem. Through both case studies, we show potential net benefits in design, open-loop scheduling, and closed-loop scheduling settings; this positions risk-constrained scheduling as a powerful option to abate risk in energy storage systems.

Despite the evident benefits of the risk-constrained optimization approach deployed herein, limitations remain in the scaling of computational effort and number of discretized price scenarios. This may be alleviated through the use of surrogates as done by \citet{alcantara2025quantile}, who used a quantile neural network to approximate second-stage expected value and CVaR in stochastic optimization. Using this approach can achieve a better balance between parsimony and distributional fidelity. Extending the proposed CVaR-constrained framework to endogenous price formation, e.g., through bi-level or equilibrium-based formulations, is another important direction for future work. This could be coupled with approaches for estimating more complex skewed energy prices as in \citet{matsumoto2022pricing} such that the effect of the energy price prior on risk is better understood. Lastly, we assumed participation in the DA and ID market herein; ancillary markets such as frequency containment reserves and automatic frequency restoration reserves could be jointly considered in this formulation along with uncertainties in renewable generation.

\section*{Acknowledgments}
The authors gratefully acknowledge funding from the bp International Centre for Advanced Materials (ICAM). LMPG and CT also acknowledge  support from the EPSRC (grant EP/X025292/1). CT was supported by a BASF/Royal Academy of Engineering Senior Research Fellowship. 


\nocite{*}   
\bibliographystyle{elsarticle-num-names}
\bibliography{citations} 


\pagebreak
\makenomenclature


\renewcommand\nomgroup[1]{%
  \item[\bfseries
  \ifstrequal{#1}{A}{}{%
  \ifstrequal{#1}{B}{Abbreviations}{%
  \ifstrequal{#1}{C}{Greek variables}{%
  \ifstrequal{#1}{D}{Latin variables}{%
  \ifstrequal{#1}{E}{Script variables}{%
  \ifstrequal{#1}{F}{Subscripts}{%
  \ifstrequal{#1}{G}{Superscripts}{}}}}}}}%
]}

\newcommand{\nomunit}[1]{%
\renewcommand{\nomentryend}{\hspace*{\fill}#1}}
\renewcommand{\nompreamble}{\begin{multicols}{2}}
\renewcommand{\nompostamble}{\end{multicols}}

\mbox{}

\nomenclature[B, 01]{BESS}{Battery energy storage system}
\nomenclature[B, 02]{(C)VaR}{(Conditional) Value-at-Risk}
\nomenclature[B, 03]{DA}{Day-ahead}
\nomenclature[B, 04]{DRI}{Direct reduced iron}
\nomenclature[B, 05]{EVPI}{Expected value of perfect information}
\nomenclature[B, 06]{EEV}{Expected solution of the expected value problem}
\nomenclature[B, 07]{IHS}{Integrated hydrogen system}
\nomenclature[B, 08]{ID}{Intraday}
\nomenclature[B, 09]{ISO-NEW}{Independent systems operator - New England}
\nomenclature[B, 10]{mn}{Million}
\nomenclature[B, 10]{k}{Thousand}
\nomenclature[B, 10]{SAA}{Sample average approximation}
\nomenclature[B, 11]{SOC}{State of charge}
\nomenclature[B, 12]{SOH}{State of health}
\nomenclature[B, 13]{SP}{Stochastic problem}
\nomenclature[B, 14]{VSS}{Value of stochastic solution}
\nomenclature[B, 15]{WS}{Wait and see}

\nomenclature[C, 01]{$\alpha$}{Risk percentile}
\nomenclature[C, 02]{$\epsilon$}{Risk tolerance}
\nomenclature[C, 03]{$\zeta$}{Value-at-risk}
\nomenclature[C, 04]{$\eta$}{Auxiliary variable}
\nomenclature[C, 05]{$\xi$}{Uncertainty}
\nomenclature[C, 06]{$\pi$}{Probability}
\nomenclature[C, 07]{$\sigma$}{Standard deviation}
\nomenclature[C, 08]{$\Psi$}{Support function}

\nomenclature[D, 01]{$\textbf{c}/c$}{First-stage costs/charge decisions}
\nomenclature[D, 02]{$d$}{Dispatch decisions (IHS)/discharge decisions (BESS)}
\nomenclature[D, 03]{\textbf{f}}{System model}
\nomenclature[D, 04]{$n$}{Number of}
\nomenclature[D, 05]{$p$}{Pressure}
\nomenclature[D, 06]{$t$}{time}
\nomenclature[D, 07]{$v$}{Discretized second-stage cost function}

\nomenclature[D, 08]{$C$}{Unit capacities}
\nomenclature[D, 09]{$F$}{Flowrate}
\nomenclature[D, 10]{$I$}{Inventory (IHS)/binary variable (BESS)}
\nomenclature[D, 11]{$J$}{Objective function term}
\nomenclature[D, 12]{$P$}{Prices}
\nomenclature[D, 13]{$V$}{Second-stage cost function}
\nomenclature[D, 14]{\textbf{X}}{First-stage decisions}
\nomenclature[D, 16]{\textbf{Y}}{Second-stage decisions}

\nomenclature[E, 01]{$\mathcal{C}$}{Capacity (IHS)/charge (BESS) decision feasible region}
\nomenclature[E, 02]{$\mathcal{D}$}{Dispatch (IHS) /discharge (BESS) decision feasible region}
\nomenclature[E, 03]{$\mathcal{F}$}{Cumulative distribution function}
\nomenclature[E, 04]{$\mathcal{L}$}{Objective function}
\nomenclature[E, 05]{$\mathcal{M}$}{Electricity market set}
\nomenclature[E, 06]{$\mathcal{P}$}{Probability distribution}
\nomenclature[E, 07]{$\mathcal{S}$}{Scenario set}
\nomenclature[E, 08]{$\mathcal{T}$}{Time set}
\nomenclature[E, 09]{$\mathcal{U}$}{Unit set}
\nomenclature[E, 10]{$\mathcal{X}$}{First-stage decision feasible region}
\nomenclature[E, 11]{$\mathcal{Y}$}{Second-stage decision feasible region}

\nomenclature[F, 01]{0}{First time instance}
\nomenclature[F, 02]{cap}{Capital}
\nomenclature[F, 03]{comp}{Compressor}
\nomenclature[F, 04]{elec}{Electrolyzer}
\nomenclature[F, 05]{$f$}{Final time}
\nomenclature[F, 06]{$fc$}{Fuel cell}
\nomenclature[F, 07]{heat}{Heater}
\nomenclature[F, 08]{$m$}{Market}
\nomenclature[F, 09]{nom}{Nominal value}
\nomenclature[F, 10]{obs}{Observation time}
\nomenclature[F, 11]{$s$}{Scenario}
\nomenclature[F, 12]{stor}{Storage}

\nomenclature[G, 01]{$AC$}{Alternating current}
\nomenclature[G, 02]{$DC$}{Direct current}
\nomenclature[G, 03]{$H_2$}{Hydrogen}
\nomenclature[G, 04]{in}{Into unit}
\nomenclature[G, 05]{op}{Operating}
\nomenclature[G, 06]{out}{Out of unit}

\printnomenclature

\end{document}